\documentclass[a4paper,10pt]{scrartcl}
\pdfoutput=1
\usepackage[utf8]{inputenc}
\usepackage[T1]{fontenc}
\usepackage[english]{babel}
\usepackage{amsmath}
\usepackage{mathtools}
\usepackage{amssymb}
\usepackage{xcolor}
\allowdisplaybreaks	
\usepackage{hyperref}
\usepackage{microtype}
\usepackage{enumerate}

\newtheorem{theorem}{Theorem}
\newtheorem{corollary}{Corollary}[theorem]
\newtheorem{lemma}[theorem]{Lemma}
\newtheorem{fact}[theorem]{Fact}
\newtheorem{definition}{Definition}
\newtheorem{remark}{Remark}
\newtheorem{example}{Example}

\newcommand{\Mt}[2]{\stackrel{t}{M}\!^{#1}_{#2}}
\newcommand{\Mv}[2]{\stackrel{v}{M}\!^{#1}_{#2}}
\newcommand{\Mw}[2]{\stackrel{w}{M}\!^{#1}_{#2}}
\newcommand{\grad}[0]{\operatorname{grad}}

\newcommand{\pdiff}[2][]{\frac{\partial#1}{\partial#2}}
\newcommand{\opd}{\operatorname{d}}
\newcommand{\doi}[1]{\href{http://www.doi.org/#1}{{\footnotesize\uppercase{doi}}: \texttt{#1}}}
\newcommand{\issn}{{\footnotesize\uppercase{issn}}:}
\newcommand{\isbn}{{\footnotesize\uppercase{isbn}}:}

\title{Geodesic flows of c-projectively equivalent metrics are quantum integrable}
\author{Jan Schumm\thanks{\texttt{jan.schumm@uni-jena.de}, supported by DFG (MA 2565/4)}\\
Friedrich Schiller University, Jena}

\begin{document}
	
	\maketitle
\begin{abstract}
	Given two c-projectively equivalent metrics on a K\"ahler manifold we show that the canoncially constructed, Poisson-commuting integrals of motion of the geodesic flow, linear and quadratic in momenta, also commute as quantum operators. The methods employed here also provide a proof of a similar statement in the case of projective equivalence. 
	We also investigate the addition of potentials, i.e.\ the generalization to natural Hamiltonian systems.
	We show that the commuting operators lead to separation of variables for Schr\"odinger's equation.
\end{abstract}

\section{C-projective geometry, integrals and quantization rules}\label{cprogeo}
\begin{definition}[K\"ahler manifold]
	A K\"ahler manifold (of arbitrary signature) is a manifold $\mathcal{M}^{2n}$ of real dimension $2n$ endowed with the following objects:
	\begin{itemize}
		\item a (pseudo-)riemannian metric $g$ and its associated Levi-Civita connection $\nabla$
		\item a complex structure $J$, i.e. an endomorphism on the space of vector fields with $J^2=-Id$ 
		\item $g$ and $J$ must be compatible in the sense that $g(JX, Y)=-g(X,JY)$ and $\nabla J =0$
		\item We denote by $\Omega$ the two-form $\Omega(X,Y)=g(JX,Y)$
	\end{itemize}
\end{definition}
\begin{definition}[J-planar curves]
	A regular curve $\gamma: I\rightarrow \mathcal{M}$ is called \emph{$J$-planar} if there exist functions $A,B: I\rightarrow \mathbb{R}$ such that
	\begin{equation}
	\nabla_{\dot{\gamma}}\dot{\gamma}= A(t) \dot{\gamma} + B(t) J(\dot{\gamma})
	\end{equation}
	is fulfilled on $I$. Here $\dot{\gamma}$ denotes the tangent vector to $\gamma$. 
\end{definition}
This is a natural generalization of geodesics on (pseudo)-Riemannian manifolds that in arbitrary parametrization are solutions of the equation $\nabla_{\dot{\gamma}}\dot{\gamma}= A(t) \dot{\gamma}$. Similarly, the property of a curve to be $J$-planar survives under reparametrization.

\begin{definition}[C-projective equivalence]
	Let $g,\tilde{g}$ be two  K\"ahler metrics (of arbitrary signature) on $(\mathcal{M},J)$. They are called \emph{c-projectively equivalent} if and only if every $J$-planar curve of $g$ is also a $J$-planar curve of $\tilde{g}$. 
	(If every $J$-planar curve of $g$ is also a $\tilde{J}$-planar curve of $\tilde{g}$  and vice versa then the complex structures $J, \tilde{J}$ coincide up to a sign, so we did not restrict ourselves in defining c-projective equivalence for the case where both metrics are K\"ahler w.r.t.\ the same complex structure.)
\end{definition}
This is the K\"ahler analogue of \emph{projective equivalence} on (pseudo-)riemannian manifolds and was proposed by T.\ Otsuki and Y.\ Tashiro \cite{otsuki1954curves}. For a thorough introduction to c-projective geometry see \cite{cprojectivegeometry}. For completeness, we recall the definition of \emph{projective equivalence}, since the theorem \ref{projectivetheorem} is the projective analogue of theorem \ref{Newtheorem}. Because their proofs run in parallel all statements about the projective setting will be phrased as remarks placed after their c-projective counterparts.
\begin{definition}[Projective equivalence]
	Let $g,\tilde{g}$ be two (pseudo-)riemannian metrics on a manifold $\mathcal{M}$. They are called \emph{projectively equivalent} if and only if every unparametrized geodesic of $g$ is also an unparametrized geodesic of $\tilde{g}$.
\end{definition}

\subsection{The tensor $A$}	
\begin{fact}\cite{domashevmikes}, see also \cite[\S 5]{cprojectivegeometry}.
	Two (pseudo-)riemannian metrics  $g, \tilde{g}$ that are K\"ahler on a manifold $(\mathcal{M},J)$ are c-projectively equivalent if and only if the tensor
	\begin{equation*}
	A^i_j \stackrel{\operatorname{def}}{=}\left|\frac{\det \tilde{g}}{\det g}\right|^{\frac{1}{2(n+1)}}\tilde{g}^{il}g_{lj},\qquad \text{where} \qquad \tilde{g}^{il} \tilde{g}_{lm}=\delta^i_m
	\end{equation*} 
	satisfies the equation
	\begin{align}\label{eq: A differential formula}
	\nabla_k A_{ij} & = \lambda_i g_{jk} + \lambda_j g_{ik} +\bar{\lambda}_i \Omega_{jk} + \bar{\lambda}_j \Omega_{ik}\\
	\intertext{where}
	\lambda&\stackrel{\operatorname{def}}{=} \frac{1}{4}\operatorname{tr}A, \qquad
	\lambda_i \stackrel{\operatorname{def}}{=} \nabla_i \lambda \qquad \text{and} \qquad \bar{\lambda}_i=J^j_i \lambda_j \notag
	\end{align}
\end{fact}
Here and throughout the rest of the paper we use the Einstein sum convention. An index preceded by a comma is meant to indicate a covariant derivative. Raising and lowering indices is always by means of $g$: $\lambda^i=g^{ij}\lambda_j$, where $g^{is} g_{sj}= \delta^i_j$. A covariant \emph{(c-)projectively equivalent metric} $\tilde{g}$, as stated above will be the sole exception: $\tilde{g}^{is} \tilde{g}_{sj}=\delta^i_j$.
\begin{remark}[for the projective case]\cite{matveevbolsinovbenenti}
	Two (pseudo-)riemannian metrics  $g, \tilde{g}$ on a manifold $\mathcal{M}$ are projectively equivalent if and only if the tensor
	\begin{equation}\label{eq: A definition projective}
	A^i_j \stackrel{\operatorname{def}}{=}\left|\frac{\det \tilde{g}}{\det g}\right|^{\frac{1}{n+1}}\tilde{g}^{il}g_{lj},\qquad \text{where} \qquad \tilde{g}^{il} \tilde{g}_{lm}=\delta^i_m
	\end{equation} 
	satisfies the equation
	\begin{equation}\label{eq: A differential formula projective}
	\nabla_k A_{ij}= \lambda_i g_{jk} + \lambda_j g_{ik} \qquad \text{with}\qquad
	\lambda \stackrel{\operatorname{def}}{=} \frac{1}{2}\operatorname{tr}A, \qquad
	\lambda_i \stackrel{\operatorname{def}}{=} \nabla_i \lambda
	\end{equation}
\end{remark}
\begin{definition}
	We shall call Hermitian ($g$-self-adjoint and $J$-commuting) solutions of (\ref{eq: A differential formula}) \emph{c-projectively compatible} or simply \emph{c-compatible} with $(g,J)$. Likewise, symmetric solutions of (\ref{eq: A differential formula projective}) shall be called \emph{projectively compatible} or simply \emph{compatible} with $g$.
\end{definition}
\subsection{Conserved quantities of the geodesic flow}	
Throughout this paper we shall canonically identify symmetric covariant tensors with polynomials on $T^*\mathcal{M}$ via the isomorphisms $^\flat$ and $ ^\sharp$:
\begin{equation*}
\begin{split}
	^\flat :& T^{a_1\ldots a_l} \mapsto T^{a_1\ldots a_l} p_{a_1} \ldots p_{a_l}\\
	^\sharp :& T^{a_1\ldots a_l} p_{a_1} \ldots p_{a_l} \mapsto T^{(a_1\ldots a_l)} 
\end{split}
\end{equation*}
By the parentheses we mean symmetrisation with the appropriate combinatorial factor: $T^{(a_1\ldots a_l)}=1/l! \sum_{(b_1,\ldots b_l)=\pi(a_1\ldots a_l)}^{}T^{b_1 \ldots b_l}$

\noindent Let $(g,J,A)$ be c-compatible on $\mathcal{M}$ and consider the one-parameter family
\begin{equation}\label{Killingtensordefinition}
\stackrel{t}{K}\!^{ij} \stackrel{\operatorname{def}}{=}\sqrt{\det (t Id - A)}\ {(tId-A)^{-1}}^i_l g^{lj}
\end{equation}
Throughout this paper the root is to be taken in such a way that we simply halve the powers of the eigenvalues in $\det (t Id -A)$. This is well defined because all eigenvalues of $A$ are of even multiplicity. In particular $\sqrt{\det (t Id -A)}$ can be negative and it is smooth also near points where $\det (t Id -A)=0$. With the tensors $\stackrel{t}{K}$ we associate the functions $\stackrel{t}{I}: T^*\! \mathcal{M}\rightarrow \mathbb{R}$
\begin{equation}\label{eq:quadraticintegrals}
\stackrel{t}{I}\ \stackrel{\operatorname{def}}{=}\stackrel{t}{K}\!^{ij} p_i p_j
\end{equation}
\begin{theorem}P. Topalov\cite{topalov2003geodesic}, see also \cite[\S 5]{cprojectivegeometry} \label{Topalovtheorem}.
	Let $(g,J,A)$ be c-compatible. Then for any pair of real numbers $(v,w)$ the quantities $\stackrel{v}{I}$ and $\stackrel{w}{I}$ poisson-commute, i.e.\ $\{\stackrel{v}{I},\stackrel{w}{I}\}=0$\\
\end{theorem}
\begin{remark}
	The quantities $\stackrel{t}{K}$ are well-defined for all values of $t$: if we denote by $2n$ the dimension of the manifold, then $\stackrel{t}{K}$ is a polynomial of degree $n-1$. This is a consequence of $J^2=-Id$, the antisymmetry of $J$ with respect to $g$, the commutativity of $A$ with $J$ and the construction of $\stackrel{t}{K}$ from $A$.
\end{remark}
\begin{remark}[for the projective case]
	Let $(g,A)$ be compatible. We shall define ${\stackrel{t}{K}\!^{ij} \stackrel{\operatorname{def}}{=}\det (t Id - A)\ {(tId-A)^{-1}}^i_l g^{lj}}$ for the projective case. Then for any pair $s,t \in \mathbb{R}$ the quantities $\stackrel{s}{I}\ \stackrel{\operatorname{def}}{=}\stackrel{s}{K}\!^{ij} p_i p_j$ and $\stackrel{t}{I}\ \stackrel{\operatorname{def}}{=}\stackrel{t}{K}\!^{ij} p_i p_j$ are commuting integrals of the geodesic flow for $g$ \cite{matveevbolsinovbenenti, matveev1998trajectory}. 
	The projective and the c-projective case differ merely by the power of the determinant.
\end{remark}
For a c-compatible structure $(g,J,A)$ there also exists a one-parameter family of commuting integrals of the geodesic flow that are linear in momenta given by 
\begin{equation}\label{Killingvectorfielddefinition}
	\stackrel{t}{L}=\stackrel{t}{V}\!^j p_j ,\qquad \stackrel{t}{V}\!^j = J^j_k g^{ki} \nabla_i \sqrt{\det(t Id -A)}
\end{equation}
The commutation relations 
\begin{equation}
	\{\stackrel{t}{L},\stackrel{s}{I}\}=0
\end{equation}
also hold.
The number of functionally independent integrals within the family $\stackrel{t}{I}$ is equal to the degree of the minimal polynomial of $A$. The number of functionally independent integrals within $\stackrel{t}{L}$ is equal to the number of nonconstant eigenvalues of $A$. Furthermore the integrals $\stackrel{t}{I}$ are functionally independent from $\stackrel{t}{L}$, see \cite[\S 5]{cprojectivegeometry}. Thus the canonically constructed integrals are sufficient in number for Liouville-integrability of the geodesic flow if all eigenvalues of $A$ are non-constant (then they are automatically pairwise different and each eigenvalue is of multiplicity two \cite[lemma 5.16]{cprojectivegeometry}).
\subsection{Quantization rules and commutators of operators}
We adopt the quantization rules introduced by B.\ Carter \cite{carter1977} and C.\ Duval and G.\ Valent \cite[\S 3]{duvalvalent2005}, see these and their references for a reasoning and more details.
It is sufficient for our purposes to recall the quantization formulae they give: for a homogeneous polynomial $P_m: T^*\mathcal{M}\rightarrow \mathbb{R}$ of degree $m$ we construct its symmetric contravariant tensor via $^\sharp$ and  compose with the covariant derivative:
\begin{equation}\label{eq: Quantizationrules}
\begin{split}
P_0 &\mapsto \hat{P}_0 \stackrel{\operatorname{def}}{=} P_0^\sharp Id\\
P_1&\mapsto \hat{P}_1 \stackrel{\operatorname{def}}{=} \frac{i}{2}((P_1^\sharp)^j\ \circ \nabla_j + \nabla_j \circ (P_1^\sharp)^j)\\
P_2 &\mapsto \hat{P}_2 \stackrel{\operatorname{def}}{=} - \nabla_j \circ (P_2^\sharp)^{jk} \circ \nabla_k\\
P_3 &\mapsto \hat{P}_3 \stackrel{\operatorname{def}}{=} -\frac{i}{2} (\nabla_j \circ (P_3^\sharp)^{jkl} \circ \nabla_k \circ \nabla_l +\nabla_j \circ \nabla_k \circ (P_3^\sharp)^{jkl} \circ \nabla_l)
\end{split}
\end{equation}
For polynomials that are not homogeneous the quantization shall be done by quantizing the homogeneous parts and adding the results. 
So far we have been considering polynomials of degree two on the cotangent bundle of degree at most two and covariant tensors of valence at most $(2,0)$. These correspond to differential operators of degree at most 2.
But the commutator of two such second order operators generally is an operator of order three. 
Later on we can facilitate the expression for the commutator of the quantum operators of two polynomials of degree two by using the quantum operator of the poisson bracket of the two polynomials of degree two. 

\section{Results}\label{results}
The main result of this paper is a quantum version of theorem \ref{Topalovtheorem}: using the quantization rules \eqref{eq: Quantizationrules} from \cite{carter1977} and \cite{duvalvalent2005} we construct differential operators from symmetric covariant tensors 
and show that these differential operators commute.

\noindent Let $(g,J,A)$ be c-compatible and $\stackrel{t}{K}, \stackrel{t}{I}$ denote the associated Killing tensors and integrals of the geodesic flow. By \eqref{eq: Quantizationrules} their associated quantum operators are:
\begin{equation*}
\stackrel{t}{\hat{I}}(f) \stackrel{\operatorname{def}}{=} -\nabla_j \circ \stackrel{t}{K}\!^{jk} \circ \nabla_k f
\end{equation*}
(The differing letters $I$ and $K$ must not confuse the reader, for $\stackrel{t}{I}\!^\sharp = \stackrel{t}{K}$.)

\begin{theorem}\label{Newtheorem}
	Let $(g,J,A)$ be \emph{c-compatible}. Then for any pair $(v,w)$ the operators $\stackrel{v}{\hat{I}}$ and $\stackrel{w}{\hat{I}}$ commute, i.e.\ $ {[\stackrel{v}{\hat{I}},\stackrel{w}{\hat{I}}]=0}$\\
\end{theorem}	
This is a new result for both the riemannian as well as the pseudo-riemmanian case.	
	
	\begin{theorem}\label{projectivetheorem}
		Let $(g,A)$ be \emph{projectively} compatible. Then for any pair $(v,w)$ the operators $\stackrel{v}{\hat{I}}$ and $\stackrel{w}{\hat{I}}$ commute, i.e.\ $ {[\stackrel{v}{\hat{I}},\stackrel{w}{\hat{I}}]=0}$
	\end{theorem}
\begin{remark} Theorem \ref{projectivetheorem} was already proven by V. Matveev  \cite{matveev2019quantum, matveev2001quantum}. The proof that we give however only uses $C^3$-smoothness whereas the original proof used $C^8$. The proof that will be given here runs in parallel with the proof of theorem \ref{Newtheorem}. A series of remarks to the proof of theorem \ref{Newtheorem} will thus provide the proof of theorem \ref{projectivetheorem}, giving the intermediate steps for the projective case and pointing out the analogues and differences. 
\end{remark}
	
\noindent We then improve the result of theorem \ref{Newtheorem} by adding potential terms to these second order differential operators, finding commuting quantum observables for certain natural hamiltonian systems.
\begin{theorem}\label{KommutierendePotentiale}
	Let $(g,J,A)$ be c-compatible. Let 
	\begin{equation*}
	\stackrel{t}{\hat{I}}\ \stackrel{\operatorname{def}}{=} - \nabla_j\ \circ \stackrel{t}{K}\!^{jk} \circ \nabla_k ,  \qquad 
	\stackrel{t}{K}\!^{ij} \stackrel{\operatorname{def}}{=}\sqrt{\det (t Id - A)}\ {(tId-A)^{-1}}^i_l g^{lj}
	\end{equation*} be as in theorem \ref{Newtheorem}. 
	
	\noindent Let 
	\begin{equation}
	\stackrel{t}{\hat{L}}=\frac{i}{2} ( \nabla_j \circ \stackrel{t}{V}\!^j + \stackrel{t}{V}\!^j \circ \nabla_j), \qquad  \stackrel{t}{V}\!^j = J^j_k g^{ki} \nabla_i \sqrt{\det(t Id -A)}
	\end{equation}
	
	\noindent be the differential operators associated with the canonical Killing vector fields of $g$.
	Let $\stackrel{\operatorname{nc}}{E}=\{\varrho_1, \ldots, \varrho_r\}$ be the set of non-constant eigenvalues of $A$. Let $\stackrel{\operatorname{c}}{E}=\{\varrho_{r+1}, \ldots, \varrho_{r+R}\}$ be the set of constant eigenvalues and $E=\stackrel{\operatorname{nc}}{E}\cup\stackrel{\operatorname{c}}{E}$. Denote by $m(\varrho_{i})$ the algebraic multiplicity of $\varrho_i$. Let the family of potentials $\stackrel{t}{U}$, parametrized by $t$, be given by
	
	\begin{equation}\label{potential}
	\stackrel{t}{U}=\sum_{i=1}^{r+R} \prod_{\varrho_l \in E \setminus\{\varrho_i\}}^{} \frac{(t-\varrho_l)^{m(\varrho_l)/2}}{(\varrho_i-\varrho_l)^{m(\varrho_l)/2}} (t-\varrho_i)^{m(\varrho_i)/2-1} f_i
	\end{equation}
		with $\opd f_l \circ A = \varrho_l \opd f_l$ for all $l=1\ldots r+R$ and with $\opd f_l$ proportional to $\opd \varrho_l$ for all $l$ for which $\varrho_l$ is non-constant. Let associated operators $\stackrel{t}{\hat{U}}$ act on functions by mere multiplication, i.e.\ for any point $p$ we have $(\stackrel{t}{\hat{U}}(f))(p)=\stackrel{t}{U}(p) f(p)$.
	Then the operators 
	\begin{equation}
	\stackrel{t}{\hat{Q}}\ \stackrel{\operatorname{def}}{=}\ \stackrel{t}{\hat{I}}+\stackrel{t}{\hat{U}}, \qquad \stackrel{t}{\hat{L}}
	\end{equation}
	commute within the one-parameter-families as well as crosswise, i.e. for all values of $t,s \in \mathbb{R}$:
	\begin{equation}
		[\stackrel{t}{\hat{Q}},\stackrel{s}{\hat{Q}}] = [\stackrel{t}{\hat{Q}},\stackrel{s}{\hat{L}}] = [\stackrel{t}{\hat{L}},\stackrel{s}{\hat{L}}] = 0
	\end{equation}
\end{theorem}
\begin{remark}\label{killingvectorssimplifydiffoperators}
	Since $\stackrel{t}{V}\!^j = J^j_k g^{ki} \nabla_i \sqrt{\det(t Id -A)}$ is a Killing vector field for any choice of the real parameter $t$ \cite[\S 5]{cprojectivegeometry}, we have $\stackrel{t}{\hat{L}}=\frac{i}{2} ( \nabla_j \circ \stackrel{t}{V}\!^j + \stackrel{t}{V}\!^j \circ \nabla_j)=i \stackrel{t}{V}\!^j \circ \nabla_j$.
\end{remark}
\begin{remark}\ 
	\begin{enumerate}
		\item We do not discuss whether the $\stackrel{t}{U}$ are smooth at all points of the manifold. Smoothness is guaranteed at points that have a neighbourhood in which the number of different eigenvalues is constant (see the definition \ref{regularpointdef} of regular points below), provided that the $f_i$ are smooth.
		\item Formula \eqref{potential} generally allows $\stackrel{t}{U}$ to be complex-valued. 
		The conditions under which $\stackrel{t}{U}$ is real for any choice of $t\in \mathbb{R}$ are the following: for any real eigenvalue $\varrho_i$ of $A$ the corresponding function $f_i$ must be real-valued. For all pairs $(\varrho_i, \varrho_j=\bar{\varrho_i})$ of complex-conjugate eigenvalues of $A$ the corresponding functions $f_i$ and $f_j$ must be complex conjugate to each other: $f_i = \bar{f}_j$.
		\item The potentials that are admissible to be added to the quantum operators are the same that may be added to the Poisson commuting integrals. In the proof we show that the quantization imposes no stronger conditions on the potential than classical integrability and then use the the conditions imposed by the poisson-brackets to find the allowed potentials.
		
	\end{enumerate}
\end{remark}
\begin{theorem}\label{PotentialForm}
	Let $(g,J,A)$ be c-compatible and $A$ semi-simple. Let $\stackrel{t}{\hat{I}}$, $\stackrel{t}{\hat{L}}$, be as in theorem \ref{KommutierendePotentiale}.
	Then, for the operators 
	\begin{equation}
	\stackrel{t}{\hat{Q}}\ \stackrel{\operatorname{def}}{=}\ \stackrel{t}{\hat{I}}+\stackrel{t}{\hat{U}}, \qquad \stackrel{t}{\hat{L}}
	\end{equation}
	the commutation relations $
	[\stackrel{t}{\hat{Q}},\stackrel{s}{\hat{Q}}] = [\stackrel{t}{\hat{Q}},\stackrel{s}{\hat{L}}] = 0$ are satisfied if and only if the potentials are of the form (\ref{potential}) with the sole exception that a function of $t$ alone may be added to $\stackrel{t}{U}$.
\end{theorem}
\begin{corollary}
	Let $(g,J,A)$, $\stackrel{t}{\hat{I}}$, $\stackrel{t}{\hat{L}}$, $\stackrel{t}{\hat{U}}$ be as in theorem \ref{KommutierendePotentiale}. Let $\hat{I}_{(l)}, \hat{L}_{(l)}, \hat{U}_{(l)}$ be the coefficients of $t^l$ in $\stackrel{t}{\hat{I}}, \stackrel{t}{\hat{L}}, \stackrel{t}{\hat{U}}$ respectively. Then the commutation relations
	\begin{multline}\label{mixedommutators}
	[\stackrel{t}{\hat{I}}+\stackrel{t}{\hat{U}}, \hat{I}_{(l)} + \hat {U}_{(l)}] = [\stackrel{t}{\hat{I}}+\stackrel{t}{\hat{U}}, \hat{L}_{(l)}] = [\stackrel{t}{\hat{L}}, \hat{I}_{(m)} + \hat {U}_{(m)}]\\*
	=  [\stackrel{t}{\hat{L}}, \hat{L}_{(l)}]= [\stackrel{t}{\hat{L}}, \hat{I}_{(m)} + \hat {U}_{(m)}] = 0
	\end{multline}
	\begin{equation}\label{equivalentcommutators}
	[\hat{I}_{(l)} + \hat {U}_{(l)}, \hat{I}_{(m)} + \hat {U}_{(m)}] = [\hat{L}_{(l)},\hat{L}_{(m)}] = [\hat{I}_{(l)} + \hat{U}_{(l)}, \hat{L}_{(m)}] = 0
	\end{equation}
	hold true for any value of $t$ and any values $l,m \in \{1, \ldots, n-1\}$. \\
	Equations (\ref{equivalentcommutators}) are equivalent to 
	\begin{equation}\label{equivalentcommutators2}
	[\stackrel{t}{\hat{I}} + \stackrel{t}{\hat{U}}, \stackrel{s}{\hat{I}} + \stackrel{s}{\hat{U}}] = [\stackrel{t}{\hat{L}},\stackrel{s}{\hat{L}}] = [\stackrel{t}{\hat{I}} + \stackrel{t}{\hat{U}}, \stackrel{s}{\hat{L}}] = 0
	\end{equation}
	Equations (\ref{mixedommutators}), (\ref{equivalentcommutators}), (\ref{equivalentcommutators2}) remain true if a function of $t$ alone is added to $\stackrel{t}{\hat{U}}$ and constants $c_{(l)}$ are added to $\hat{U}_{(l)}$.
	If all eigenvalues of $A$ are non-constant and $\stackrel{t}{\hat{I}}, \stackrel{t}{\hat{V}}$ are as in theorem \ref{KommutierendePotentiale} then no other than the described $\hat{U}_{(l)}$ can be found such that the commutation relations above hold.
\end{corollary}
Lastly, we shall show how the search for common eigenfunctions of the operators can be reduced to differential equations in lower dimension in appropriate coordinates around regular points. In particular if all eigenvalues of $A$ are non-constant we can reduce it to ordinary differential equations only. Moreover: the case where all eigenvalues of $A$ are non-constant provide an example of $\emph{reduced separability}$ of Schr\"odinger's equation as described in \cite{BenentiChanuRastelli2002}.
\begin{definition}\label{regularpointdef}
	Let $(g,J,A)$ be c-compatible. A point $x \in \mathcal{M}$ is called regular with respect to $A$ if in a neighbourhood of $x$ the number of different eigenvalues of $A$ is constant and for each eigenvalue $\varrho$ either $\opd \varrho \neq 0$ or $\varrho$ is constant in a neighbourhood of $x$. The set of regular points shall be denoted $\mathcal{M}^0$.
\end{definition}
\begin{definition}
	Let $(g,J,A)$ be c-compatible on $\mathcal{M}$. 
	A local normal coordinate system for $\mathcal{M}$ is a coordinate system where $(g,J,A)$ assume the form of example \ref{example:general}. Existence of such coordinates in the neighbourhood of regular points is guaranteed by theorem 1.6 in \cite{bolsinov2015localc}.
\end{definition}
	
\begin{example}\label{example:general}[General example for c-compatible structures $(g,J,\omega,A)$]\cite[Example 5]{bolsinov2015localc}
	Let $2n\geq 4$ and consider an open subset $W$ of $\mathbb{R}^{2n}$ of the form $W=U\times V \times S_1 \times \ldots\times S_L \times S_{L+1} \times\ldots \times S_{L+Q}$ for open subsets $V,U\subseteq \mathbb{R}^r$, $S_\gamma \subseteq \mathbb{R}^{4m_{c_\gamma}}$ for $\gamma=1,\ldots, L$ and $S_\gamma \subseteq \mathbb{R}^{2m_{c_\gamma}}$ for $\gamma=L+1, \ldots, L+Q$.
	Let the coordinates on $U$ be separated into $l$ complex coordinates $z_1 \ldots z_l$ and $q$ real coordinates $x_{l+1}, \ldots x_{l+q}$ and introduce the tuple $(\chi_1, \ldots \chi_r)=(z_1, \bar{z}_1, \ldots, z_l, \bar{z}_l, x_{l+1}, \ldots x_{l+q})$. 
	Suppose the following data is given on these open subsets
	\begin{itemize}
		\item K\"ahler structures $(g_\gamma, J_\gamma, \omega_\gamma)$ on $S_\gamma$ for $\gamma=1, \ldots, L+Q$
		\item For each $\gamma=1, \ldots, L+Q$, a parallel hermitian endomorphism $A_\gamma : TS_\gamma \rightarrow TS_\gamma$ for $(g_\gamma, J_\gamma)$. For $\gamma=1\ldots L$ $A$ has a pair of complex conjugate eigenvalues $c_\gamma, \bar{c}_\gamma$ of equal algebraic multiplicity $m(c_\gamma)=m(\bar{c}_\gamma)$. For $\gamma=L+1\ldots L+Q$ $A$ has a single real eigenvalue $c_\gamma$ of algebraic multiplicity $m(c_\gamma)$. 
		\item  Holomorphic functions $\sigma_j(z_j)$ for $1\leq j\leq l$ and smooth functions $\sigma_j(x_j)$ for $l+1\leq j\leq r$.
	\end{itemize}
	Moreover, we choose 1-forms $\alpha_1, \ldots \alpha_r$ on $S=S_1\times\cdots \times S_N$ that satisfy
	\begin{equation}
		\opd \alpha_i = (-1)^i \sum_{\gamma=1}^{L+Q} \omega_{\gamma} (A^{r-i}_\gamma \cdot, \cdot)
	\end{equation}
	To facilitate the expressions for the c-compatible structure that will be constructed, the following expressions shall be introduced:\,the tuple $E\!=\!(\varrho_1,\ldots, \varrho_n)\!=\!(\sigma_1, \bar{\sigma}_1,\ldots, \sigma_l, \bar{\sigma}_l,\allowbreak \sigma_{l+1}, \ldots, \sigma_{l+q}, c_1, \bar{c_1}, \ldots ,\allowbreak c_L,\allowbreak \bar{c}_L,\allowbreak c_{L+1}, \ldots, c_{L+Q})$ contains the designated eigenvalues for $A$. Their algebraic multiplicities shall be denoted by $(m(\varrho_l),l\!=\!1,\ldots , r+R)\!=\!(2,\ldots,2 , m(c_1), m(\bar{c}_1), \ldots)$. The non-constant eigenvalues shall be collected in order in $\stackrel{\operatorname{nc}}{E}=(\varrho_1,\ldots, \varrho_r)=(\sigma_1, \bar{\sigma_1},\ldots,\allowbreak \sigma_l, \ldots \bar{\sigma}_l, \sigma_l+1, \ldots, \sigma_{l+q})$ and the collection of constant eigenvalues shall be referenced as $\stackrel{\operatorname{c}}{E}=E\setminus\stackrel{\operatorname{nc}}{E}$.
		The quantity $\Delta_i$ for $i=1\ldots r$ is given by $\Delta_i = \prod_{\varrho \in \stackrel{\operatorname{nc}}{E}\setminus \{\varrho_i\}}^{}(\varrho_i -\varrho)$. The function $\mu_i$ denotes the elementary symmetric polynomial of degree $i$ in the variables $\stackrel{\operatorname{nc}}{E}$ and $\mu_i(\hat{\varrho_s})$ denotes the elementary symmetric polynomial of degree $i$ in the variables $\stackrel{\operatorname{nc}}{E}\setminus \{\varrho_i\}$. We shall further define the one-forms $\vartheta_1, \ldots \vartheta_r$ on $W$ via $\vartheta_i = \opd t_i +\alpha_i$.\medskip
		
	\noindent Suppose that at every point of $W$ the elements of $\stackrel{\operatorname{nc}}{E}$ are mutually different and different from the constants $c_1, \bar{c}_1, \ldots , c_R$ and their differentials are non-zero. Then $(g,\omega, J)$ given by the formulae

	\begin{equation}\label{eq:example general 1}
	\begin{split}
	g=& \sum_{i=1}^{r} \varepsilon_i \Delta_i \opd \chi_i^2 + \sum_{i=0}^r (-1)^i \mu_i \sum_{\gamma=1}^{L+Q} g_\gamma (A_\gamma^{r-i} \cdot, \cdot)\\
	&+ \sum_{i,j=1}^l \left[ \sum_{s=1}^r \frac{\mu_{i-1} (\hat{\varrho}_s) \mu_{j-1}(\hat{\varrho}_s)}{\varepsilon_s \Delta_s} \left(\pdiff[\varrho_s]{\chi_s}\right)^2\right]\vartheta_i \vartheta_j\\
	\omega =& \sum_{i=1}^r  \opd \mu_i \wedge \vartheta_i + \sum_{i=0}^{r} (-1)^i \mu_i \sum_{\gamma=1}^{L+Q} \omega_\gamma ( A_\gamma^{r-i} \cdot, \cdot)\\
	\end{split}
	\end{equation}
	\begin{equation}\label{eq:example general 2}
	\begin{split}
		\opd \chi_i \circ J =& -\frac{1}{\varepsilon_i \Delta_i} \pdiff[\varrho_i]{\chi_i} \sum_{j=1}^r \mu_{j-1} (\hat{\varrho_i})\vartheta_i, \qquad 
		\vartheta_i \circ J = (-1)^{i-1} \sum_{j=1}^r \varepsilon_j \varrho_j^{r-i} \left(\pdiff[\varrho_j]{\chi_j}\right)^{-1} \opd \chi_j
	\end{split}
	\end{equation}	
	is K\"ahler, where $(\varepsilon_1, \ldots \varepsilon_{2l}, \varepsilon_{2l+1}, \ldots \varepsilon_r)=(-1/4, \ldots, -1/4, \pm 1, \ldots, \pm 1)$ determine the signature of $g$.
	
	\noindent With local coordinates $\stackrel{\gamma}{y}$ on $S_\gamma$ we write $\alpha_i = \sum_{\gamma,q}^{} \stackrel{\gamma}{\alpha}_{iq} \opd \stackrel{\gamma}{y}_q$ and\\ ${A_\gamma= \sum_{p,q}  (A_\gamma)^q_p \opd \stackrel{\gamma}{y}_p \otimes \partial_{\stackrel{\gamma}{y}_q}}$. Then the endomorphism $A$ given by	
	\begin{multline}\label{eq:example general 3}
		A= \sum_{s=1}^{r} \varrho_s \opd \chi_s \otimes \partial_{\chi_s} + \sum_{i,j=1}^{r} (\mu_i \delta_{1j} - \delta_{i(j-1)}) \vartheta_i \otimes \partial_{t_j} \\
		+\sum_{\gamma=1}^{L+Q}\sum_{p,q}^{} (A_\gamma)^q_p \opd \stackrel{\gamma}{y}_p \otimes \left(\partial_{\stackrel{\gamma}{y}_q}- \sum_{i=1}^{r} \stackrel{\gamma}{\alpha}_{iq} \partial_{t_i}\right)
	\end{multline}
	is c-compatible with $(g,J,\omega)$.
\end{example}

\begin{theorem}\cite[Theorem 1.6 / Example 5]{bolsinov2015localc}
	Suppose (g,J,A) are c-compatible on $\mathcal{M}$ of real dimension $2n$. \\
	Assume that in a
	small neighbourhood $W \subseteq \mathcal{M}^0$ of a regular point, A has
	\begin{itemize}
		\item $r=2l+q$ non-constant eigenvalues on $W$ which separate into $l$ pairs of complex-conjugate
		eigenvalues $\varrho_1 ,\bar{\varrho}_1 , \ldots , \varrho_{l} , \bar{\varrho}_{l} : W \rightarrow \mathbb{C}$ and $q$ real eigenvalues $\varrho_{l+1} , \ldots ,\allowbreak \varrho_{l+q} : W \rightarrow \mathbb{R}$,
		\item $R=2L+Q$ constant eigenvalues which separate into $L$ pairs of complex conjugated eigenvalues $c_1, \bar{c}_1, \ldots , c_L,  \bar{c}_L$ and $Q$ real eigenvalues $c_{L+1}, \ldots , c_{L+Q}$
	\end{itemize}
	then the K\"ahler structure $(g,J,\omega)$ and $A$ are given on $W$ by the formulas of example \ref{example:general}. 
\end{theorem}

\begin{theorem}\label{eigenfunctionstheorem}
	Let $(g,J,A)$ be c-compatible on $\mathcal{M}$. Let $A$ be semi-simple and let all constant eigenvalues be real.
	Let $(g,J,\omega,A)$ be given by the formulas of example \ref{example:general} and adopt the naming conventions of example \ref{example:general}.

	Let $\psi$ be a simultaneous eigenfunction of 
	\begin{equation*}
	\stackrel{t}{\hat{Q}}\ \stackrel{\operatorname{def}}{=} - \nabla_j\ \circ \stackrel{t}{K}\!^{jk} \circ \nabla_k + \stackrel{t}{\hat{U}},  \qquad 
	\stackrel{t}{K}\!^{ij} \stackrel{\operatorname{def}}{=}\sqrt{\det (t Id - A)}\ {(tId-A)^{-1}}^i_l g^{lj}
	\end{equation*}
	and
	\begin{equation*}
	\stackrel{t}{\hat{L}}=\frac{i}{2} ( \nabla_j \circ \stackrel{t}{V}\!^j + \stackrel{t}{V}\!^j \circ \nabla_j), \qquad  \stackrel{t}{V}\!^j = J^j_k g^{ki} \nabla_i \sqrt{\det(t Id -A)}
	\end{equation*}
	for all $t$, where
	\begin{equation}
	\stackrel{t}{U}=\sum_{i=1}^{r} \prod_{l=1, l \neq i}^{r} \frac{(t-\varrho_l)^{m(\varrho_l)/2}}{(\varrho_i-\varrho_l)^{m(\varrho_l)/2}} (t-\varrho_i)^{m(\varrho_i)/2-1} f_i
	\end{equation}
	with $\opd f_i \circ A = \varrho_i \opd f_i$. Then there exist constants $\tilde{\lambda}_0, \ldots, \tilde{\lambda}_{r+R-1}, \omega_1,\ldots,\omega_r$, such that $\psi$ satisfies the following ODE:
	\begin{equation}\label{eq:separated ode}
	\begin{split}
	\frac{-1}{\varepsilon_k \varrho_k '}\partial_{\chi_k} \varrho_k' \prod_{\varrho_{\gamma} \in \stackrel{\operatorname{c}}{E}} (\varrho_{k} - \varrho_{\gamma})^{m(\varrho_{\gamma})/2} \partial_{\chi_k} \psi \hspace{5cm}  & \\
	 + \sum_{i,j=1}^{r} \frac{\varepsilon_k (-\varrho_k)^{2r-i-j}}{(\varrho_k')^2}\prod_{\varrho_{\gamma}\in \stackrel{\operatorname{c}}{E}}^{}(\varrho_k - \varrho_{\gamma})^{m(\varrho_{\gamma})/2} \omega_i \omega_j
	+ f_r \psi = & \sum_{i=0}^{n-1} \lambda_i \varrho_k^i\\
	i \partial_{t_q} \psi = & \omega_q \psi
	\end{split}
	\end{equation}
	for $q=1\ldots r$, where the $\lambda_i$ are given by $\sum_{i=0}^{n-1} \lambda_i s^i=\prod_{\varrho_{\gamma}\in \stackrel{\operatorname{c}}{E}}^{} (s-\varrho_{\gamma})^{m(\varrho_{\gamma})/2-1}\sum_{j=0} ^{r+R-1}\tilde{\lambda}_j$ and $\psi$ also fulfills the partial differential equations:
	\begin{equation}\label{eq:separated pde}
	\begin{split}
	\sum_{j=0}^{r+R-1} \tilde{\lambda}_j \varrho_{\gamma}^j
	=&
		-\,\prod_{\mathclap{\varrho_c \in \stackrel{\operatorname{c}}{E}\setminus \{\varrho_\gamma\}}}^{} (\varrho_\gamma - \varrho_c) 
		\left[ \frac{1}{|\det g_\gamma |^{1/2}} \partial_{\stackrel{\gamma}{y_i}} g_\gamma^{ij} |\det g_\gamma |^{1/2} \partial_{\stackrel{\gamma}{y_j}} \psi\right.\\
		&\left. -i \sum_{q=1}^{r} \frac{1}{|\det g_\gamma |^{1/2}} \partial_{\stackrel{\gamma}{y_i}} g_\gamma^{ij} |\det g_\gamma |^{1/2} \stackrel{\gamma}{\alpha}_{qj} \omega_q \psi -i \sum_{q=1}^{r} g_\gamma^{ij} \stackrel{\gamma}{\alpha}_{qi} \omega_q \partial_{\stackrel{\gamma}{y_j}} \psi\right.\\
		&\left. - \sum_{p,q=1}^{r} g_\gamma^{ij} \stackrel{\gamma}{\alpha}_{qi}\stackrel{\gamma}{\alpha}_{pj} \omega_q \omega_p \psi\right]
		+ \frac{1}{\prod_{\varrho_c \in \stackrel{\operatorname{c}}{E}\setminus\{\varrho_\gamma\}}(\varrho_\gamma - \varrho_c)^{m(\varrho_c)/2-1}} f_\gamma \psi
	\end{split}
	\end{equation}
	for $\gamma=r+1\ldots r+R$. \\
	The converse is also true: If a function $\psi$ satisfies equations \eqref{eq:separated ode} and \eqref{eq:separated pde} for some constants $\tilde{\lambda}_0, \ldots, \tilde{\lambda}_{r+R-1}, \omega_1,\ldots,\omega_r$, then it is an eigenfunction of $\stackrel{t}{\hat{Q}}$ and $\stackrel{t}{\hat{L}}$.
\end{theorem}
A particular application of theorem \ref{eigenfunctionstheorem} is in the treatment of the Laplace-Beltrami operator in the case where the number of integrals is maximal: if $\mathcal{M} $ is compact and riemannian its eigenfunctions provide means to construct a countable basis in $L^2(\mathcal{M})$, see e.g. \cite{jostgeometricanalysis}. Because the integrals are self-adjoint and do commute pairwise as well as with the Laplace-Beltrami-operator they leave each other's eigenspaces invariant and thus  there exists a countable basis of $L^2(\mathcal{M})$ that consists of simultanous eigenfunctions of the commuting operators. In our case, since $\Delta$ is within the span of our one-parameter family of operators we can find such basis by solving the simultaneous eigenvalue problems corresponding to the integrals only. In the case where all eigenvalues of $A$ are non-constant this reduces to ODE.

\section{Proof of the results}\label{proofs}

\subsection{Basic facts}\label{basicfacts}
We shall provide some formulae that will be used throughout the proof of theorem \ref{Newtheorem}. We will be mainly working with (\ref{eq: A differential formula}). Unless it is stated otherwise we are working in the \emph{c-projective} setting with only some remarks providing the analogue formulae for the \emph{projective setting}.

We start by computing the covariant derivative of $\det A$. By using Jacobi's formula $\opd \det M= \operatorname{tr}(\operatorname{Ad}(M) \opd M)$, we get:
\begin{equation}\label{eq: ddetA}
\begin{split}
	\nabla_k \det (A) =& \det (A)  {A^{-1}}^q_r \nabla_k  A^r_q\\
	=& \det (A) {A^{-1}}^q_r \left[ \lambda^r g_{qk} + \lambda_q \delta^r_k + g^{rs} \bar{\lambda}_s \Omega_{qk} + \bar{\lambda}_q g^{rs} \Omega_{sk} \right] \\
	=& 4 \det (A) {A^{-1}}^s_k \lambda_s
\end{split}
\end{equation}
Equation (\ref{eq: A differential formula}) was used to expand the covariant derivative of $A$ and then the properties of $J$ interacting with $g$ were exploited to obtain this.
\begin{lemma}(\cite[\S 5]{cprojectivegeometry})\label{lemma: A selfadjoint lambda}
	$A$ and thus $t Id - A$ and $(t Id-A)^{-1}$ are self-adjoint with respect to $\nabla^2 \lambda$: ${A^{-1}}^j_l\lambda_{j,k} = {A^{-1}}^j_k\lambda_{j,l}$. The last claim is of course only true if $(t Id-A)^{-1}$ exists.
\end{lemma}
That is, if we consider the second derivative of $\lambda$ in the way that it maps two vector fields $\xi, \eta$ to a scalar function on $\mathcal{M}$ then
we have $\nabla^2 \lambda(A\xi,\eta)=\nabla^2 \lambda(\xi,A\eta)\forall (\xi,\eta)$.\\
\textbf{Proof of lemma \ref{lemma: A selfadjoint lambda}}: Without loss of generality we may assume that $A$ is invertible. Otherwise we may choose $\varepsilon$ such that $\varepsilon Id - A$ is invertible. Then we can apply the same procedure that we will apply to $\ln \det A$ to $\ln \det (\varepsilon Id -A)$ instead. We can therefrom show that $(\varepsilon Id -A)^{-1}$ is self adjoint with respect to $\nabla^2 \lambda$. And thus by linear algebra $\varepsilon Id -A$ and consequentially $A$ are also $\nabla^2 \lambda$-self-adjoint.\\
Now to the main argument: We compute the second covariant derivative of $\ln \det A$ using (\ref{eq: ddetA}), (\ref{eq: A differential formula}), the general identity $\opd A^{-1}=A^{-1}\cdot (\opd  A) \cdot A^{-1}$ as well as the antisymmetry of $J$ with respect to $g$:
\begin{equation}
\begin{split}
\frac{1}{4} \nabla_k \nabla_l \ln\det A =& \lambda_j {A^{-1}}^j_p g^{ps} \lambda_s {A^{-1}}^q_l g_{qk} + \lambda_j  {A^{-1}}^j_k {A^{-1}}^q_l \lambda_q \\ 
&+ \lambda_j {A^{-1}}^j_p g^{ps} \bar{\lambda}_s {A^{-1}}^q _l \Omega_{qk} - \lambda_j {A^{-1}}^j_p J^p_k \bar{\lambda}_q {A^{-1}}^q_l
+{A^{-1}}^j_l \lambda_{j,k}
\end{split}
\end{equation}
The left hand side is symmetric with respect to $(k,l)$. The first, second and fourth term on the right hand side are symmetric as well. The third term vanishes. Consequently the last term must be symmetric as well. Thus $A^{-1}$ is self-adjoint with respect to $\nabla^2 \lambda$.
By means of linear algebra, the self-adjointness with respect to $\lambda_{j,k}$ is also true for $(t Id -A)$ or $(t Id-A)^{-1}$. The latter of course is only true, provided that $t$ is not chosen to be within the spectrum of $A$. This concludes the proof of lemma \ref{lemma: A selfadjoint lambda}.
\begin{remark}[for the projective setting]
	The formula $A^j_l\lambda_{j,k} = A^j_k\lambda_{j,l}$ is true for the projective setting as well. Performing the same computations as in lemma \ref{lemma: A selfadjoint lambda} gives the intermediate results $\nabla_k \det (t Id - A)
	= - 2 \det (t Id -A) {(t Id - A)^{-1}}^s_k \lambda_s$ and $\frac{1}{2}\nabla_k \nabla_l \ln \det A = \lambda_j {A^{-1}}^j_p g^{ps} \lambda_s {A^{-1}}^q_l g_{qk} + \lambda_j  {A^{-1}}^j_k {A^{-1}}^q_l \lambda_q +{A^{-1}}^j_l \lambda_{j,k}$ to which the same logic is applied as in the c-projective case.
\end{remark}
\begin{lemma}\label{lemma: R A S selfadjoint}
	Let $S$ be an endomorphism on the space of vector fields on $\mathcal{M}$ with the following properties: $J\circ S=S\circ J$, $A\circ S = S \circ A$,\ $\nabla^2 \lambda(S\cdot,\cdot)=\nabla^2 \lambda(\cdot,S\cdot)$ and $g(S\cdot,\cdot)=g(\cdot,S\cdot)$. Then the formula
	\begin{equation}\label{eq: R A S selfadjoint}
	R^r_{\ ijk} A_{rl} S^{ij} - R^r_{\ ijl} A_{rk} S^{ij}= 0
	\end{equation}
	is valid.
\end{lemma}
The sign of the Riemann tensor is chosen such that $R^i_{\ jkl} = \partial_k \Gamma^i_{lj} -\partial_l \Gamma^i_{kj} +\Gamma^i_{ks}\Gamma^s_{lj} - \Gamma^i_{ls}\Gamma^s_{kj}$. It is obviously not essential for the formula but the proof would need a few changes in the signs.\\

\noindent\textbf{Proof of lemma \ref{lemma: R A S selfadjoint}}
In the proof we reuse and extend the ideas used in the proof of  equation (12) in \cite{kiosakmatveev2010}.
We inspect the second derivative of $A$:
\begin{equation}\label{eq: A second derivative expanded}
\begin{split}
R^r_{\ jkl} A_{ir} + R^r_{\ ikl} A_{rj} = & (\nabla_l \nabla_k - \nabla_k \nabla_l) A_{ij} \\
	= & \lambda_{i,l} g_{jk} - \lambda_{i,k} g_{jl} +\lambda_{j,l} g_{ik} -\lambda_{j,k}g_{il}\\
	& +\bar{\lambda}_{i,l} \Omega_{jk} - \bar{\lambda}_{i,k} \Omega_{jl} +\bar{\lambda}_{j,l} \Omega_{ik} - \bar{\lambda}_{j,k}\Omega_{il}
\end{split}
\end{equation}
The first equality is the Ricci identity and is true for any $(0,2)$-tensor. The second equality comes from differentiating (\ref{eq: A differential formula}).
We continue by adding the equation with itself three times after performing cyclic permutations of $(j,k,l)$. The three terms rising from the first term on the left hand side of \eqref{eq: A second derivative expanded} vanish due to the Bianchi identity. On the right hand side only terms involving $\bar{\lambda}$ remain:
\begin{equation}\label{eq: A second derivative cyclic}
\begin{split}
	R^r_{\ ikl} A_{rj} + R^r_{\ ijk} A_{rl} + R^r_{\ ilj} A_{rk}= & (
	+\bar{\lambda}_{i,l} \Omega_{jk} + \bar{\lambda}_{i,k} \Omega_{lj} + \bar{\lambda}_{i,j} \Omega_{kl} \\
	&- \bar{\lambda}_{i,k} \Omega_{jl} - \bar{\lambda}_{i,l} \Omega_{kj} - \bar{\lambda}_{i,j} \Omega_{lk}) + (i\leftrightarrow j)
\end{split}
\end{equation}
$+(i \leftrightarrow j)$ means that the preceding term is to be added again, but with indices $i$ and $j$ interchanged. We now multiply this equation with $S^{ij}$. Since $S$ is $g$-self-adjoint and commutes with $J$ we have $\bar{\lambda}_{i,j} S^{ij}=0$. Using this and using again that $S$ commutes with $J$ the right hand side simplifies as follows:
\begin{equation*}
	\begin{split}
	R^r_{\ ikl} A_{rj} S^{ij} + R^r_{\ ijk} A_{rl} S^{ij}+ R^r_{\ ilj} A_{rk} S^{ij} = 4 (S^j_k \lambda_{j,l} - S^j_l \lambda_{j,k})
	\end{split}
\end{equation*}
$S\circ A$ is g-self-adjoint because $S$ and $A$ are g-self-adjoint and commute. As a consequence, we have that $R^r_{\ ikl} A_{rj} S^{ij}$ vanishes on the left hand side. Further utilizing the symmetry of the curvature tensor as well as the self-adjointness of $S$ with respect to $\nabla^2 \lambda$ we reach the desired result:
\begin{equation*}
R^r_{\ ijk} A_{rl} S^{ij} - R^r_{\ ijl} A_{rk} S^{ij}= 0
\end{equation*}
\begin{remark}[for the projective case]
	In the projective case lemma \ref{lemma: R A S selfadjoint} also holds true. The procedure involves the same steps as in the proof of lemma \ref{lemma: R A S selfadjoint}. Removing the terms involving $\Omega$ from \eqref{eq: A second derivative expanded} gives the intermediate step for the projective case. Equation \ref{eq: A second derivative cyclic} trivially simplifies to $R^r_{\ ikl} A_{rj} + R^r_{\ ijk} A_{rl} + R^r_{\ ilj} A_{rk}=0$. After multiplication with $S^{ij}$ the first term vanishes with the same argument as in the c-projective case and the result is obtained. 
	Since the proof of corollary \ref{corollary: Riemann Ricci selbstadjungiert} only involves linear algebra the arguments are exactly the same in the projective and the c-projective case.
\end{remark}
\begin{corollary}\label{corollary: Riemann Ricci selbstadjungiert}Linear algebra applied to (\ref{eq: R A S selfadjoint}) provides the formula:
\begin{equation}\label{eq: Riemann Ricci selbstadjungiert}
	{(v Id - A)^{-1}}^{lr} R^k_{\ jir}S^{ij} - {(v Id - A)^{-1}}^{kr} R^{l}_{\ jir}S^{ij}=0
\end{equation}

\end{corollary}
\textbf{Proof}:	Consider the endomorphism $ R^r_{\ ijl} S^{ij}$ on the space of vector fields. Raising and lowering indices in (\ref{eq: R A S selfadjoint}) shows that it commutes with $A$. Consequently it also commutes with $(v Id -A)$ and thus with $(v Id -A)^{-1}$. Standard index manipulations using the symmetries of the curvature tensor imply the result.

\subsection{Proof of theorem \ref{Newtheorem}}
The family $\stackrel{t}{\hat{I}}$ is a polynomial in $t$ and therefore continuous. It is therefore sufficient to show that the commutator vanishes for all $v$ and $w$ that are are not in the spectrum of $A$. Otherwise we can consider two sequences $(v_n)_{n\in \mathbb{N}}$ and $(w_n)_{n\in \mathbb{N}}$ where none of the elements of the sequence are in the spectrum of $A$ and that converge to $v$ and $w$. Then for each of the pairs $(v_n,w_m)$ from the sequences the commutator $[\stackrel{v_n}{\hat{I}},\stackrel{w_m}{\hat{I}}]$ will vanish and consequently it will vanish in the limit $(m,n)\rightarrow\infty$.\\
Equations (3.12), (3.13) and (3.14) in \cite{duvalvalent2005} give us the general formula for the commutator of two operators formed from arbitrary homogeneous polynomials $P_2, Q_2$ of degree two on $T^*\mathcal{M}$:
\begin{equation}\label{Commutator}
\begin{split}
\left[\hat{P}_2 , \hat{Q}_2\right]=&i\widehat{\{P_2,Q_2\}}+\frac{2}{3} (\nabla_j B^{jk}_{P_2 , Q_2})\nabla_l\\
\end{split}
\end{equation}
Here $\{P_2,Q_2\}$ is the Poisson bracket of the two homogeneous polynomials $P_2$ and $Q_2$ (of degree two) on $T^*\mathcal{M}$. $\{P_2,Q_2\}$ is a polynomial of degree three in momenta. We recall that the ``hat'' over $\{P_2,Q_2\}$ is explained in \eqref{eq: Quantizationrules}:
This polynomial is mapped to a differential operator according to
\begin{equation*}
	\hat{\cdot}: P_3 \mapsto \hat{P}_3 \stackrel{\operatorname{def}}{=} -\frac{i}{2} (\nabla_j \circ P_3^{jkl} \circ \nabla_k \circ \nabla_l +\nabla_j \circ \nabla_k \circ P_3^{jkl} \circ \nabla_l)
\end{equation*}
where for a given polynomial  $P_3$ the quantities $P_3^{jkl}$ are chosen such that they are symmetric and $P_3= P_3^{jkl}p_j p_k p_l$.

The tensor $B^{kl}_{P_2 , Q_2}$ is given by the formula 
\begin{multline}\label{eq: B allgemein}
B^{jk}_{P_2 , Q_2}=P^{l\left[j\right.}\nabla_l \nabla_m Q^{\left. k\right] m} 
+ P^{l\left[j\right.} R^{\left. k\right]}_{\ mnl} Q^{mn} 
- (P \leftrightarrow Q) \\
- \nabla_l P^{m\left[j\right.}\nabla_m Q^{\left.k\right] l}
-P^{l\left[j\right.}R_{lm} Q^{\left. k\right] m} 
\end{multline}
The brackets around the indices mean taking the antisymmetric part. The subtraction of $(P\leftrightarrow Q)$ is meant to act upon the two leftmost terms. For the two rightmost terms the antisymmetrization w.r.t.\ $(j \leftrightarrow k)$ is the same as if one were to antisymmetrize these terms w.r.t.\ $(P\leftrightarrow Q)$. For formula \eqref{eq: B allgemein} the sign $R^i_{\ jkl} = \partial_k \Gamma^i_{lj} -\partial_l \Gamma^i_{kj} +\Gamma^i_{ks}\Gamma^s_{lj} - \Gamma^i_{ls}\Gamma^s_{kj}$ is important. But the reader may forget about it at once because it is not needed for our further investigations, as will be seen in the upcoming lemmas \ref{Kruemmungstensortermverschwindet} and \ref{Riccitensorverschwindet}.

\noindent We plug the operators $\stackrel{v}{\hat{I}}$ and $\stackrel{w}{\hat{I}}$ into formula (\ref{Commutator}). Using theorem \ref{Topalovtheorem} by Topalov we get:
\begin{equation}
[\stackrel{v}{\hat{I}},\stackrel{w}{\hat{I}}]=\frac{2}{3} (\nabla_j B^{jk}_{\stackrel{v}{I} , \stackrel{w}{I}})\nabla_k
\end{equation}
\begin{remark}$i\widehat{\{P_2,Q_2\}}$ in formula (\ref{Commutator}) is a differential operator of order $3$ while the other term on the right hand side is a differential operator of order one. Therefore it is a necessary condition for the quantities $\stackrel{v}{I}$ and $\stackrel{w}{I}$ to poisson-commute in order for their associated differential operators to commute. This is of course a long known fact.
\end{remark}
\begin{remark}[for the projective case]
	The fact that $\stackrel{t}{K}$ is a family of Killing tensors polynomial of degree ${n-1}$ in $t$ and that their corresponding quadratic polynomials on $T^*\mathcal{M}$ Poisson commute pairwise can be found in \cite{matveevbolsinovbenenti}. Employing this instead of theorem \ref{Topalovtheorem} brings proof in the projective case to the point where only equation \eqref{BTensor} needs to be verified.
\end{remark}

\noindent It remains to prove that  
\begin{multline}\label{BTensor}
\nabla_j B^{jk}_{\stackrel{v}{I} , \stackrel{w}{I}}=\nabla_j\left((\stackrel{v}{K}\!^{l\left[j\right.}\nabla_l \nabla_m \stackrel{w}{K}\!^{\left. k\right] m} - (v\leftrightarrow w) )
- \nabla_l \stackrel{v}{K}\!^{m\left[j\right.}\nabla_m \stackrel{w}{K}\!^{\left.k\right] l}\right.\\
\left.+ (\stackrel{v}{K}\!^{l\left[j\right.} R^{\left. k\right]}_{\ mnl} \stackrel{w}{K}\!^{mn} 
- (v\leftrightarrow w) )
-\stackrel{v}{K}\!^{l\left[j\right.}R_{lm} \stackrel{w}{K}\!^{\left. k\right] m}\right)=0
\end{multline}
The proof will be split into three steps:\ first we show that $ \stackrel{v}{K}\!^{l\left[j\right.} R^{\left. k\right]}_{\ mnl}\! \stackrel{w}{K}\!^{mn} {- (v\leftrightarrow w)\!=\!0}$. In the second step we show that $\stackrel{v}{K}\!^{l\left[j\right.}R_{lm} \stackrel{w}{K}\!^{\left. k\right] m}=0$. This will be done in lemma \ref{Kruemmungstensortermverschwindet} and \ref{Riccitensorverschwindet}. These will reduce (\ref{BTensor}) to 
\begin{equation}\label{BTensorreduced}
\nabla_j\left(\stackrel{v}{K}\!^{l\left[j\right.}\nabla_l \nabla_m \stackrel{w}{K}\!^{\left. k\right] m} - (v\leftrightarrow w) 
- \nabla_l \stackrel{v}{K}\!^{m\left[j\right.}\nabla_m \stackrel{w}{K}\!^{\left.k\right] l}\right)=0
\end{equation}
which we will show in the last step.
\begin{lemma}\label{Kruemmungstensortermverschwindet}
	$ \stackrel{v}{K}\!^{l\left[j\right.} R^{\left. k\right]}_{\ mnl} \stackrel{w}{K}\!^{mn} - (v\leftrightarrow w) =0$
\end{lemma}
\textbf{Proof of lemma \ref{Kruemmungstensortermverschwindet}:} It follows from lemma \ref{lemma: A selfadjoint lambda} that in corollary \ref{corollary: Riemann Ricci selbstadjungiert} we may take ${S=\sqrt{\det (w Id -A)}} {(w Id -A)^{-1}}$, i.e.\ $S^{ij}=\stackrel{w}{K}\!^{ij}$. Plugging this into equation (\ref{eq: Riemann Ricci selbstadjungiert}) and multiplying with $\sqrt{\det(v Id - A)}$ gives:
\begin{equation}\label{kruemmungstensorverschwindetfast}
	{2 \stackrel{v}{K}\!^{l\left[j\right.} R^{\left. k\right]}_{\ mnl} \stackrel{w}{K}\!^{mn}=0}
\end{equation}
Interchanging $v$ and $w$, subtracting the result from this and dividing by 2 proves lemma \ref{Kruemmungstensortermverschwindet}.
\begin{remark}[for the projective case]
	Lemma 13 is true in the projective and the c-projective case. The proof for the projective case only requires to take \\${S=\det (w Id -A)} {(wId-A)^{-1}}$ and multiplying with $\det (v Id -A)$ instead of\\ $\sqrt{\det (v Id -A)}$.
\end{remark}

\begin{lemma}\label{Riccitensorverschwindet}
	$\stackrel{v}{K}\!^{l\left[j\right.}R_{lm}\! \stackrel{w}{K}\!^{\left. k\right] m}=0$
\end{lemma}
\textbf{Proof of lemma \ref{Riccitensorverschwindet}}
If we let $S=Id$ in corollary \ref{corollary: Riemann Ricci selbstadjungiert}, then in formula (\ref{eq: Riemann Ricci selbstadjungiert}) the multiplication with $S^{ij}=g^{ij}$ means contraction of the Riemann tensor to the negative of the Ricci tensor. Raising and lowering indices yields that $(v Id - A)^{-1}$ commutes with the Ricci tensor when both are considered as endomorphisms on the space of vector fields:
\begin{equation}\label{AkommutiertmitRicci}
		{(v Id - A)^{-1}}\!^{r}_l R^k_{r} - {(v Id - A)^{-1}}\!^{k}_r R^{r}_{l}=0
\end{equation}
Of course $(w Id -A)^{-1}$ commutes with $(v Id - A)^{-1}$ and the Ricci tensor as well, so by multiplying (\ref{AkommutiertmitRicci}) with ${(w Id - A)^{-1}}\!^l_j$ and using the commutativity gives
\begin{equation}\label{eq:36}
	{(w Id - A)^{-1}}\!^k_r R^r_{l} {(v Id - A)^{-1}}\!^{l}_j - {(v Id - A)^{-1}}\!^{k}_r R^{r}_{l} {(w Id - A)^{-1}}\!^l_j=0
\end{equation}
After multiplication of \eqref{eq:36} with $\sqrt{\det (v Id -A) \det (w Id -A)}$ and raising and lowering indices, lemma \ref{Riccitensorverschwindet} is proven.
\begin{remark}[for the projective case]
	\hspace{-.4 em} Replacing the multiplication of $\sqrt{{\det (v Id -\!A)}}\allowbreak \sqrt{\det (w Id -A)}$ with $\det (v Id -A) \det (w Id -A)$ is the only difference between the proof of lemma \ref{Riccitensorverschwindet} in the projective and c-projective case.
\end{remark}

\noindent Having established lemma \ref{Kruemmungstensortermverschwindet} and \ref{Riccitensorverschwindet} we now compute the terms ${\stackrel{v}{K}\!^{l\left[j\right.}\nabla_l \nabla_m \stackrel{w}{K}\!^{\left. k\right] m}}\allowbreak - (v \leftrightarrow w)$ and $\nabla_m\stackrel{v}{K}\!^{l\left[j\right.}\nabla_l  \stackrel{w}{K}\!^{\left. k\right] m}$ separately and then show that \eqref{BTensorreduced} is fulfilled to prove theorem \ref{Newtheorem}.

\noindent We introduce a shorthand notation:
\begin{equation*}
\begin{split}
\Lambda=&\operatorname{grad}\lambda=(g^{ij} \lambda_j), \qquad \bar{\Lambda}=J \operatorname{grad}\lambda, \qquad \Mt{}{}=(t Id -A)^{-1}
\end{split}
\end{equation*}
for any $t$ outside the spectrum of $A$.\\
We get
\begin{equation*}
	\nabla_k \det (t Id - A) = - 4 \det (t Id -A) \Mt{s}{k} \lambda_s
\end{equation*}
in the same way as we have obtained formula \eqref{eq: ddetA}.
Using this, the general matrix identity $\opd C^{-1} = C^{-1}\cdot \opd C \cdot C^{-1}$ and \eqref{eq: A differential formula}, the covariant derivative of the Killing-tensor $\stackrel{t}{K}$ evaluates to:
\begin{equation}\label{Killingtensorderivative}
\begin{split}
\nabla_k \stackrel{t}{K}\! ^{jl} =&\sqrt{\det(t Id -A)} [-2 \Mt{s}{k} \lambda_s\ \Mt{j}{r}g^{rl}\\
&+\Mt{j}{p} \lambda^p \Mt{l}{k}+\Mt{j}{k}\ \Mt{l}{q}\lambda^q-\Mt{j}{p} g^{ps} \bar{\lambda}_s \Mt{l}{q}J^q_k-\Mt{l}{r}g^{rq}\bar{\lambda}_q \Mt{j}{p}J^p_k]
\end{split}
\end{equation}
Contracting the indices $k$ and $l$ we get:
\begin{equation}
\begin{split}
\nabla_k \stackrel{t}{K}\! ^{jk} =&\sqrt{\det(t Id -A)} [-2 \Mt{j}{l}\Mt{l}{s} \lambda^s+\Mt{j}{l}\lambda^l \operatorname{tr}(\Mt{}{})]
\end{split}
\end{equation}
For the derivative of this expression we receive
\begin{equation}\label{Kzweiteableitung}
\begin{split}
\nabla_l\nabla_m \stackrel{t}{K}\!^{km} =&\sqrt{\det(t Id -A)} [2\Mt{r}{l}\lambda_r \lambda_s g^{ts} \Mt{m}{t} \Mt{k}{m}+2 \Mt{q}{l}\Mt{r}{q}\lambda_r \lambda_s g^{ts} \Mt{k}{t}\\
&+2\Mt{r}{l}\bar{\lambda}_r \bar{\lambda}_s g^{ts} \Mt{m}{t} \Mt{k}{m}+2 \Mt{q}{l}\Mt{r}{q}\bar{\lambda}_r \bar{\lambda}_s g^{ts} \Mt{k}{t}\\
&-\Mt{r}{l}\lambda_r\lambda_s g^{st}\Mt{k}{l} \operatorname{tr}(\Mt{}{})-\Mt{r}{l}\bar{\lambda}_r\bar{\lambda}_s g^{st}\Mt{k}{l} \operatorname{tr}(\Mt{}{})\\
&-2 \Mt{k}{m}\ g(\Mt{}{}\Lambda,\Mt{}{}\Lambda)-2 \Mt{k}{r}\Mt{r}{l}\ g(\Mt{}{}\Lambda,\Lambda)+\Mt{k}{l}\ g(\Mt{}{}\Lambda,\Lambda) \operatorname{tr}(\Mt{}{})\\
&-2 \lambda_{s,l} g^{ps} \Mt{p}{t}\Mt{k}{p}+\lambda_{s,l}g^{ps} \Mt{k}{p}\operatorname{tr}(\Mt{}{})]
\end{split}
\end{equation}
We shall denote by $(\Mw{}{}\Lambda)^k$ the $k$-th component of $\Mw{}{}(\Lambda)$. Now multiplying the previous equation with $\stackrel{v}{K}\!^{lj}$ and antisymmetrizing with respect to $(j \leftrightarrow k)$ and $(v \leftrightarrow w)$ gives
\begin{equation}\label{Bersterterm}
\begin{split}
\stackrel{v}{K}\!^{l\left[j\right.}\nabla_l \nabla_m \stackrel{w}{K}\!^{\left. k\right] m} - (v \leftrightarrow w) =&\frac{1}{2}\sqrt{\det(v Id-A)}\sqrt{\det(w Id-A)}\\
&\cdot[[[2 (\Mv{}{} \Mw{}{}\Lambda)^j (\Mw{}{}\!^2 \Lambda)^k+2(\Mv{}{}\Mw{}{}\!^2 \Lambda)^j (\Mw{}{}\Lambda)^k\\
&-(\Mv{}{}\Mw{}{}\Lambda)^j (\Mw{}{}\Lambda)^k \operatorname{tr} (\Mw{}{})]\\
&+(\Lambda\leftrightarrow\bar{\Lambda})]-(j\leftrightarrow k)] -(v \leftrightarrow w) 
\end{split}
\end{equation}
Here $(\Lambda\leftrightarrow\bar{\Lambda})$ indicates that the previous bracket shall be added with $\Lambda$ replaced by $\bar{\Lambda}$, $(j \leftrightarrow k)$ indicates antisymmetrization with respect to $j$ and $k$, likewise for $(v\leftrightarrow w)$. In (\ref{Bersterterm}) the terms from \eqref{Kzweiteableitung} involving second derivatives of $\lambda$ have cancelled out as a consequence of lemma \ref{lemma: A selfadjoint lambda}. 
When forming the right hand side expression of \eqref{Bersterterm} the terms of the second to last row of \eqref{Kzweiteableitung} cancel each other out after the antisymmetrization $(j\leftrightarrow k)$ due to $(v Id-A)$ and $(w Id-A)$ commuting and being self-adjoint to $g$.

\begin{remark}[for the projective case]
	To get the formula for $\stackrel{v}{K}\!^{l\left[j\right.}\nabla_l \nabla_m \stackrel{w}{K}\!^{\left. k\right] m} - (v \leftrightarrow w)$ in the projective case we perform the same steps, using (\ref{eq: A differential formula projective}) instead of the c-projective formula (\ref{eq: A differential formula}). The next three equations give projective analogues of the formulae (\ref{Killingtensorderivative}), (\ref{Kzweiteableitung}) and (\ref{Bersterterm}):
	\begin{align}
	\nabla_k \stackrel{t}{K}\! ^{jl} =&\det(t Id -A) [-2 \Mt{s}{k} \lambda_s\ \Mt{j}{r}g^{rl}
	+\Mt{j}{p} \lambda^p \Mt{l}{k}+\Mt{j}{k}\ \Mt{l}{q}\lambda^q] \label{eq: projective calc 1}\\
	\nabla_l\nabla_m \stackrel{t}{K}\!^{km} =&\det(t Id -A) [\Mt{r}{l}\lambda_r \lambda_s g^{ts} \Mt{m}{t} \Mt{k}{m}+ \Mt{q}{l}\Mt{r}{q}\lambda_r \lambda_s g^{ts} \Mt{k}{t}\notag\\
	&	\begin{aligned}
	&-\Mt{r}{l}\lambda_r\lambda_s g^{st}\Mt{k}{l} \operatorname{tr}(\Mt{}{})\\
	&- \Mt{k}{m}\ g(\Mt{}{}\Lambda,\Mt{}{}\Lambda)- \Mt{k}{r}\Mt{r}{l}\ g(\Mt{}{}\Lambda,\Lambda)\\
	&+\Mt{k}{l}\ g(\Mt{}{}\Lambda,\Lambda) \operatorname{tr}(\Mt{}{})\\
	&- \lambda_{s,l} g^{ps} \Mt{p}{t}\Mt{k}{p}+\lambda_{s,l}g^{ps} \Mt{k}{p}\operatorname{tr}(\Mt{}{})]\end{aligned}\\
	\stackrel{v}{K}\!^{l\left[j\right.}\nabla_l \nabla_m \stackrel{w}{K}\!^{\left. k\right] m}\phantom{=}&\notag\\
	- (v \leftrightarrow w) =&\frac{1}{2}\det(v Id-A)\det(w Id-A)\label{projectivecalc1}\\
	&\begin{aligned}
	& \cdot[[(\Mv{}{} \Mw{}{}\Lambda)^j (\Mw{}{}\!^2 \Lambda)^k + 
	(\Mv{}{}\Mw{}{}\!^2 \Lambda)^j (\Mw{}{}\Lambda)^k\\
	&-(\Mv{}{}\Mw{}{}\Lambda)^j (\Mw{}{}\Lambda)^k \operatorname{tr} (\Mw{}{})] \\
	&-(j\leftrightarrow k)] -(v \leftrightarrow w) \end{aligned}\notag
	\end{align}
\end{remark}

\noindent By means of the general matrix identity
\begin{equation*}
	(v Id -A)^{-1}\!-(w Id - A)^{-1}\!=(w-v) \cdot (v Id -A)^{-1} \cdot (w Id - A)^{-1}
\end{equation*}
 which reads $\Mv{}{} - \Mw{}{}= (w-v) \Mv{}{}\Mw{}{}$ in our shorthand, as well as the trace applied to this matrix identity we expand (\ref{Bersterterm}):
\begin{align}\label{Berstertermumgeformt}
	\allowdisplaybreaks
	\stackrel{v}{K}\!^{l\left[j\right.}\nabla_l \nabla_m \stackrel{w}{K}\!^{\left. k\right] m} - (v \leftrightarrow w)  & =   \frac{1}{2}\sqrt{\det(v Id-A)}\sqrt{(w Id-A)} \cdot[[[(w-v)^{-1} \notag\\ 
	& \begin{aligned} & \phantom{=}[2 (\Mv{}{} \Lambda)^j (\Mw{}{}\!^2 \Lambda)^k - 2 (\Mw{}{} \Lambda)^j (\Mw{}{}\!^2 \Lambda)^k \\
	& \phantom{=}+2(\Mv{}{}\Mw{}{} \Lambda)^j (\Mw{}{}\Lambda)^k - 2(\Mw{}{}\!^2 \Lambda)^j (\Mw{}{}\Lambda)^k\\
	& \phantom{=}-(\Mv{}{}\Lambda)^j (\Mw{}{}\Lambda)^k \operatorname{tr} (\Mw{}{})
	+(\Mw{}{}\Lambda)^j (\Mw{}{}\Lambda)^k \operatorname{tr} (\Mw{}{})]]\\
	& \phantom{=}-(v \leftrightarrow w) ]-(j\leftrightarrow k)] +(\Lambda\leftrightarrow\bar{\Lambda})\end{aligned}\notag\\
	\intertext{We strike out terms that cancel after antisymmetrization w.r.t.$\ (j\leftrightarrow k)$:}
	&\begin{aligned} =&\frac{1}{2}\sqrt{\det(v Id-A)}\sqrt{\det(w Id-A)} \cdot[[[ (w-v)^{-1}\\
	& \cdot [2 (\Mv{}{} \Lambda)^j (\Mw{}{}\!^2 \Lambda)^k 
	+ 2(\Mv{}{}\Mw{}{} \Lambda)^j (\Mw{}{}\Lambda)^k \\
	& -(\Mv{}{}\Lambda)^j (\Mw{}{}\Lambda)^k \operatorname{tr} (\Mw{}{})]\\
	&  -(v \leftrightarrow w) ]-(j\leftrightarrow k)] +(\Lambda\leftrightarrow\bar{\Lambda})\end{aligned}\notag
	\intertext{We expand $(v \leftrightarrow w)$ and $(j \leftrightarrow k)$. The sign of the antisymmetrization is caught in the prefactor $(w-v)^{-1}$:}
	&\begin{aligned}=&\frac{1}{2}\sqrt{\det(v Id-A)}\sqrt{\det(w Id-A)} (w-v)^{-1}\\
	&\cdot[2 (\Mv{}{} \Lambda)^j (\Mw{}{}\!^2 \Lambda)^k 
	+ 2 (\Mw{}{} \Lambda)^j (\Mv{}{}\!^2 \Lambda)^k \\
	&+ 2(\Mv{}{}\Mw{}{} \Lambda)^j (\Mw{}{}\Lambda)^k 
	+ 2(\Mw{}{}\Mv{}{} \Lambda)^j (\Mv{}{}\Lambda)^k \\
	&- 2 (\Mv{}{} \Lambda)^k (\Mw{}{}\!^2 \Lambda)^j 
	- 2 (\Mw{}{} \Lambda)^k (\Mv{}{}\!^2 \Lambda)^j \\
	&- 2(\Mv{}{}\Mw{}{} \Lambda)^k (\Mw{}{}\Lambda)^j 
	- 2(\Mw{}{}\Mv{}{} \Lambda)^k (\Mv{}{}\Lambda)^j \\
	&-(\Mv{}{}\Lambda)^j (\Mw{}{}\Lambda)^k \operatorname{tr} (\Mw{}{}) - (\Mw{}{}\Lambda)^j (\Mv{}{}\Lambda)^k \operatorname{tr} (\Mv{}{})\\
	&+(\Mv{}{}\Lambda)^k (\Mw{}{}\Lambda)^j \operatorname{tr} (\Mw{}{}) + (\Mw{}{}\Lambda)^k (\Mv{}{}\Lambda)^j \operatorname{tr} (\Mv{}{})]\\
	& +(\Lambda\leftrightarrow\bar{\Lambda})\end{aligned}\notag\\
	\intertext{We can now apply $\Mv{}{} - \Mw{}{}= (w-v) \Mv{}{}\Mw{}{}$ and $\operatorname{tr}(\Mv{}{}) - \operatorname{tr}(\Mw{}{})= (w-v) \operatorname{tr}(\Mv{}{}\Mw{}{})$ in the opposite direction as before, pairing terms (1,8), (2,7), (3,6), (4,5), (9,12), (10,11) in the bracket:}
	& \begin{aligned}=&\frac{1}{2}\sqrt{\det(v Id-A)}\sqrt{\det(w Id-A)}\\
	& \cdot[[-2 (\Mv{}{}\!^2 \Mw{}{}\Lambda)^j (\Mw{}{} \Lambda)^k-2(\Mv{}{}\Lambda)^j (\Mv{}{}\Mw{}{}\!^2\Lambda)^k \\ 
	& +(\Mv{}{}\Lambda)^j (\Mw{}{}\Lambda)^k \operatorname{tr} (\Mv{}{}\cdot\Mw{}{})]\\
	&-(j\leftrightarrow k)]+(\Lambda\leftrightarrow\bar{\Lambda})\end{aligned}
\end{align}
\begin{remark}[for the projective case]
	Performing the same steps on (\ref{projectivecalc1}) gives
	\begin{equation}\label{projectivecalc2}
		\begin{split}
		\stackrel{v}{K}\!^{l\left[j\right.}\nabla_l \nabla_m \stackrel{w}{K}\!^{\left. k\right] m} - (v \leftrightarrow w)
		=&\frac{1}{2}\det(v Id-A)\det(w Id-A)\\
		&\cdot[[- (\Mv{}{}\!^2 \Mw{}{}\Lambda)^j (\Mw{}{} \Lambda)^k-(\Mv{}{}\Lambda)^j (\Mv{}{}\Mw{}{}\!^2\Lambda)^k \\ 
		& +(\Mv{}{}\Lambda)^j (\Mw{}{}\Lambda)^k \operatorname{tr} (\Mv{}{}\cdot\Mw{}{})]-(j\leftrightarrow k)]\\
		\end{split}
	\end{equation}
	for the projective scenario.
\end{remark}

\noindent We have now worked $\stackrel{v}{K}\!^{l\left[j\right.}\nabla_l \nabla_m \stackrel{w}{K}\!^{\left. k\right] m} - (v \leftrightarrow w)$ into a suitable  form. From (\ref{Killingtensorderivative}) we now compute $\nabla_m\stackrel{v}{K}\!^{l\left[j\right.}\nabla_l  \stackrel{w}{K}\!^{\left. k\right] m}$:
\begin{equation}\label{Bzweiterterm}
\begin{split}
\nabla_m\stackrel{v}{K}\!^{l\left[j\right.}\nabla_l  \stackrel{w}{K}\!^{\left. k\right] m} =&\frac{1}{2}\sqrt{\det(v Id-A)}\sqrt{\det(w Id-A)}\\
&\cdot[[-2 (\Mv{}{}\!^2 \Mw{}{}\Lambda)^j (\Mw{}{} \Lambda)^k-2(\Mv{}{}\Lambda)^j (\Mv{}{}\Mw{}{}\!^2\Lambda)^k\\
&+(\Mv{}{}\Lambda)^j (\Mw{}{}\Lambda)^k \operatorname{tr} (\Mv{}{}\cdot\Mw{}{})]\\
&-(j\leftrightarrow k)]-(\Lambda\leftrightarrow\bar{\Lambda})
\end{split}
\end{equation}
In this equation $(j\leftrightarrow k)$ yields the same result as $(v\leftrightarrow w)$.\\
\begin{remark}[for the projective case]
	By means of \eqref{eq: projective calc 1} the projective analogue of \eqref{Bzweiterterm} evaluates to
	\begin{equation}\label{projectivecalc3}
	\begin{split}
	\nabla_m\stackrel{v}{K}\!^{l\left[j\right.}\nabla_l  \stackrel{w}{K}\!^{\left. k\right] m} =&\frac{1}{2}\det(v Id-A)\det(w Id-A)\\
	&\cdot[[-2 (\Mv{}{}\!^2 \Mw{}{}\Lambda)^j (\Mw{}{} \Lambda)^k-2(\Mv{}{}\Lambda)^j (\Mv{}{}\Mw{}{}\!^2\Lambda)^k\\
	&+(\Mv{}{}\Lambda)^j (\Mw{}{}\Lambda)^k \operatorname{tr} (\Mv{}{}\cdot\Mw{}{})]-(j\leftrightarrow k)]
	\end{split}
	\end{equation}
	We see that the right hand side expression is equal to the right hand side expression of (\ref{projectivecalc2}). Thus if we plug (\ref{projectivecalc2}) and (\ref{projectivecalc3}) into (\ref{BTensorreduced}) then both terms cancel each other and (\ref{BTensorreduced}) is satisfied without even having to carry out the differentiation, concluding the proof of theorem \ref{projectivetheorem}. The fact that in the projective case $B^{jk}_{\stackrel{v}{I} , \stackrel{w}{I}}$ vanishes, whereas in the c-projective case $\nabla_j B^{jk}_{\stackrel{v}{I} , \stackrel{w}{I}}$ vanishes but $B^{jk}_{\stackrel{v}{I} , \stackrel{w}{I}}$ does not is the most significant difference between the projective and the c-projective case.
\end{remark}

We now compare (\ref{Berstertermumgeformt}) and (\ref{Bzweiterterm}): they are the same except the first is symmetric with respect to $(\Lambda \leftrightarrow \bar{\Lambda})$ while the latter is antisymmetric. Subtracting both consequently yields:
\begin{equation}
\begin{split}
\stackrel{v}{K}\!^{l\left[j\right.}\nabla_l \nabla_m \stackrel{w}{K}\!^{\left. k\right] m}  - & (v \leftrightarrow w) -\nabla_m\stackrel{v}{I}\!^{l\left[j\right.}\nabla_l  \stackrel{w}{I}\!^{\left. k\right] m}\\
 = &\sqrt{\det(v Id-A)}\sqrt{\det(w Id-A)}\\
&\cdot[-2 (\Mv{}{}\!^2 \Mw{}{}\bar{\Lambda})^j (\Mw{}{} \bar{\Lambda})^k-2(\Mv{}{}\bar{\Lambda})^j (\Mv{}{}\Mw{}{}\!^2\bar{\Lambda})^k\\ &\ +(\Mv{}{}\bar{\Lambda})^j (\Mw{}{}\bar{\Lambda})^k \operatorname{tr} (\Mv{}{}\cdot\Mw{}{})]
-(v\leftrightarrow w)
\end{split}
\end{equation}
It remains to show that \eqref{BTensorreduced} is fulfilled, that is to apply $\nabla_j$ to this expression and show that this vanishes. 
In the computation we use
\begin{itemize}
	\item the compatibility condition (\ref{eq: A differential formula})
	\item Jacobi's formula for the derivative of the determinant
	\item $\opd (A^{-1})=-A^{-1}\cdot \opd (A) \cdot A^{-1}$
\end{itemize}
to expand the left hand side expression of \eqref{BTensorreduced}. We then immediately strike out terms that vanish individually due to the self-adjointness of A with respect to $g$ and the antisymmetry of $J$ with respect to $g$:
\begin{equation}\label{halfresult}
\begin{split}
	\nabla_j(\stackrel{v}{K}\!^{l\left[j\right.}\nabla_l \nabla_m \stackrel{w}{K}\!^{\left. k\right] m} - & (v \leftrightarrow w)-\nabla_m\stackrel{v}{K}\!^{l\left[j\right.}\nabla_l  \stackrel{w}{K}\!^{\left. k\right] m}) \\
	= & -2 \sqrt{\det(v Id-A)}\sqrt{\det(w Id-A)}\\
	& \cdot [2 g(\Mv{}{}\Mw{}{}\bar{\Lambda}, \Mv{}{}\Mw{}{}\bar{\Lambda}) (\Mw{}{}\Lambda- \Mv{}{} \Lambda)^k\\
	&-2 g(\Mv{}{} \Mw{}{} \bar{\Lambda}, \Mw{}{} \bar{\Lambda} - \Mv{}{} \bar{\Lambda}) (\Mv{}{} \Mw{}{} \Lambda)^k\\
	& +2 g(\bar{\Lambda}, \Mv{}{} \Mw{}{} \bar{\Lambda})](\Mv{}{} \Mw{}{}(\Mw{}{}\Lambda- \Mv{}{} \Lambda))^k\\
	&-2 g(\bar{\Lambda}, \Mw{}{} \bar{\Lambda} - \Mv{}{} \bar{\Lambda}) (\Mv{}{}\!^2 \Mw{}{}\!^2 \Lambda)^k\\
	& + \bar{\lambda}_s g^{ls} \Mv{m}{l} \Mw{t}{m} \Mv{j}{t} \bar{\lambda}_{p,j} g^{pq} \Mw{k}{q}\\
	& - \bar{\lambda}_s g^{ls} \Mw{m}{l} \Mv{t}{m} \Mw{j}{t} \bar{\lambda}_{p,j} g^{pq} \Mv{k}{q}\\
	& + \bar{\lambda}_s g^{ps} \Mv{j}{p} \bar{\lambda}_{t,j} g^{tr} \Mw{l}{r} \Mv{m}{l} \Mw{k}{m}\\
	& - \bar{\lambda}_s g^{ps} \Mw{j}{p} \bar{\lambda}_{t,j} g^{tr} \Mv{l}{r} \Mw{m}{l} \Mv{k}{m}\\
	& + \operatorname{tr} (\Mv{}{} \Mw{}{}) \cdot[g(\bar{\Lambda}, \Mv{}{} \Mw{}{} \bar{\Lambda}) (\Mw{}{} \Lambda - \Mv{}{} \Lambda)^k\\
	& - g(\bar{\Lambda}, \Mw{}{} \bar{\Lambda} - \Mv{}{} \bar{\Lambda}) (\Mv {}{} \Mw{}{} \Lambda)^k\\
	&+ \bar{\lambda}_s g^{ps} \Mv{j}{p} \bar{\lambda}_{t,j} g^{tr} \Mw{k}{r}
	- \bar{\lambda}_s g^{ps} \Mw{j}{p} \bar{\lambda}_{t,j} g^{tr} \Mv{k}{r}]
\end{split}
\end{equation}
As a consequence of lemma \ref{lemma: A selfadjoint lambda} the terms involving second derivatives of $\lambda$ cancel each other out in this expression. The other terms cancel each other out after applying $\Mv{}{} - \Mw{}{}= (w-v) \Mv{}{}\Mw{}{}$. Thus theorem \ref{Newtheorem} is proven.\\

\subsection{Addition of potential | Proof of theorem \ref{KommutierendePotentiale} and \ref{PotentialForm}}\label{additionofpotential}
\subsubsection{Four equivalent problems}
Let $(g,J,A)$ be \emph{c-compatible}.

\begin{lemma}\label{lemma:4problems}
	Let $K= g^{ij}p_i p_j$,  \begin{equation*}
	\stackrel{t}{I}\ \stackrel{\operatorname{def}}{=} \stackrel{t}{K}\!^{jk} p_j p_k,  \qquad 
	\stackrel{t}{K}\!^{ij} \stackrel{\operatorname{def}}{=}\sqrt{\det (t Id - A)}\ {(tId-A)^{-1}}^i_l g^{lj}
	\end{equation*}
	and
	\begin{equation*}
	\stackrel{t}{L}=\stackrel{t}{V}\!^j p_j, \qquad  \stackrel{t}{V}\!^j = J^j_k g^{ki} \nabla_i \sqrt{\det(t Id -A)}
	\end{equation*}
	as well as the corresponding differential operators according to the quantization rules stated earlier.\\
	Then the following four problems are equivalent:
	describe all functions $U$, $\stackrel{t}{U}$, such that
	\begin{enumerate}
		\item\label{firstproblem} $\{\stackrel{s}{I}+\stackrel{s}{U}, K+U\}=0$ and $\{\stackrel{s}{I}+\stackrel{s}{U}, \stackrel{t}{L}\}=0 \quad \forall t,s \in \mathbb{R}$
		\item\label{secondproblem} $[\stackrel{s}{\hat{I}}+\stackrel{s}{\hat{U}}, \hat{K}+\hat{U}]=0$ and $[\stackrel{s}{\hat{I}}+\stackrel{s}{\hat{U}}, \stackrel{t}{\hat{L}}]=0 \quad \forall t,s \in \mathbb{R}$
		\item\label{thirdproblem} $\{\stackrel{s}{I}+\stackrel{s}{U}, \stackrel{t}{I}+\stackrel{t}{U}\}=0$ and $\{\stackrel{s}{I}+\stackrel{s}{U}, \stackrel{t}{L}\}=0 \quad \forall t,s \in \mathbb{R}$
		\item\label{fourthproblem} $[\stackrel{s}{\hat{I}}+\stackrel{s}{\hat{U}}, \stackrel{t}{\hat{I}}+\stackrel{t}{\hat{U}}]=0$ and $[\stackrel{s}{\hat{I}}+\stackrel{s}{\hat{U}}, \stackrel{t}{\hat{L}}]=0 \quad \forall t,s \in \mathbb{R}$
	\end{enumerate}
\end{lemma}
\textbf{Proof of lemma \ref{lemma:4problems}}:
To do so, we show that 
\begin{enumerate}[i.]
	\item\label{firstequivalence}	$\{\stackrel{t}{I}+\stackrel{t}{U}, K+U\}=0 \forall t\in \mathbb{R} \Leftrightarrow  [\stackrel{t}{\hat{I}}+\stackrel{t}{\hat{U}}, \hat{K}+\hat{U}]=0 \forall t \in \mathbb{R} \Leftrightarrow \stackrel{t}{K}\!^i_j \pdiff[U]{x^i}= \pdiff[\stackrel{t}{U}]{x^j} \forall t \in \mathbb{R}$
	
	\item\label{secondequivalence}	$\{\stackrel{t}{I}+\stackrel{t}{U}, K+U\}=0 \forall t\in \mathbb{R} \Leftrightarrow  \{\stackrel{t}{I}+\stackrel{t}{U}, \stackrel{s}{I}+\stackrel{s}{U}\}=0 \forall s,t \in \mathbb{R}$
	
	\item\label{thirdequivalence}	$\{\stackrel{t}{I}+\stackrel{t}{U}, \stackrel{s}{I}+\stackrel{s}{U}\}=0 \forall s,t \in \mathbb{R} \Leftrightarrow  [\stackrel{t}{\hat{I}}+\stackrel{t}{\hat{U}}, \stackrel{s}{\hat{I}}+\stackrel{s}{\hat{U}}]=0 \forall s,t \in \mathbb{R}$\\ $ \Leftrightarrow 
	\stackrel{t}{K}\!^i_j \pdiff[\stackrel{s}{U}]{x^i}={\stackrel{s}{K}\!^i_j} \pdiff[\stackrel{t}{U}]{x^i}  \forall s,t\in \mathbb{R}$
	
	\item\label{fourthequivalence} $\{\stackrel{t}{I}+\stackrel{t}{U}, \stackrel{s}{L}\}=0 \quad \forall t,s \in \mathbb{R} \Leftrightarrow [\stackrel{t}{\hat{I}}+\stackrel{t}{\hat{U}}, \stackrel{s}{\hat{L}}]=0 \quad \forall t,s \in \mathbb{R} \Leftrightarrow \opd \stackrel{t}{U} (\stackrel{s}{V})=0\ \forall s,t \in \mathbb{R}$
\end{enumerate}
It is implied that all equations are to hold for all choices of its parameters, we shall not specify it each and every time again.

\noindent To \ref{thirdequivalence}:  We use the linearity of the commutator:
\begin{equation}
[\stackrel{s}{\hat{I}}+\stackrel{s}{\hat{U}},\stackrel{t}{\hat{I}}+\stackrel{t}{\hat{U}}]=[\stackrel{s}{\hat{I}},\stackrel{t}{\hat{I}}]+[\stackrel{s}{\hat{I}},\stackrel{t}{\hat{U}}]+[\stackrel{s}{\hat{U}},\stackrel{t}{\hat{I}}]+[\stackrel{s}{\hat{U}},\stackrel{t}{\hat{U}}]
\end{equation}
The term $[\stackrel{s}{\hat{I}},\stackrel{t}{\hat{I}}]$ vanishes due to theorem \ref{Newtheorem} and $[\stackrel{s}{\hat{U}},\stackrel{t}{\hat{U}}]$ vanishes trivially since the operators corresponding to the potentials act merely by multiplication. 
In \cite{duvalvalent2005} or by direct computation we have $[\stackrel{s}{\hat{I}},\stackrel{t}{\hat{U}}] = \widehat{\{\stackrel{s}{I}, \stackrel{t}{U}\}}$. 
Consequently $[\stackrel{s}{\hat{I}}+\stackrel{s}{\hat{U}},\stackrel{t}{\hat{I}}+\stackrel{t}{\hat{U}}]= \widehat{\{\stackrel{s}{I}, \stackrel{t}{U}\}}+\widehat{\{\stackrel{s}{U}, \stackrel{t}{I}\}}$. Since quantization is a linear map and only the zero polynomial is mapped to a vanishing differential operator, we have that $[\stackrel{s}{\hat{Q}},\stackrel{t}{\hat{Q}}]=0$ if and only if $\{\stackrel{s}{I}, \stackrel{t}{U}\}+\{\stackrel{s}{U}, \stackrel{t}{I}\}=0$. This in turn is true if and only if $\{\stackrel{s}{I}, \stackrel{t}{U}\}^\sharp+\{\stackrel{s}{U}, \stackrel{t}{I}\}^\sharp=0$. Expressing this in terms of $\stackrel{t}{K}, \stackrel{t}{U}, \stackrel{s}{K}, \stackrel{s}{U}$ and lowering an index and rearranging terms yields $\stackrel{t}{K}\!^i_j \pdiff[\stackrel{s}{U}]{x^i}=\stackrel{s}{K}\!^i_j \pdiff[\stackrel{t}{U}]{x^i}$. Likewise using the fact that $\{\stackrel{s}{I}, \stackrel{t}{I}\}=0$, we have $\{\stackrel{t}{I}+\stackrel{t}{U}, \stackrel{s}{I}+\stackrel{s}{U}\}= \{\stackrel{s}{I}, \stackrel{t}{U}\}+\{\stackrel{s}{U}, \stackrel{t}{I}\}$.

\noindent Statement \ref{firstequivalence} can be seen analogously to \ref{thirdequivalence} since $K$ lies in the span of $\stackrel{t}{K}$.

\noindent To \ref{secondequivalence}: It suffices to show the equivalence of the rightmost equations of items \ref{firstequivalence} and \ref{thirdequivalence}:
Fix an arbitrary value for $t$ in $\stackrel{t}{K}\!^i_j \pdiff[\stackrel{s}{U}]{x^i}=\stackrel{s}{K}\!^i_j \pdiff[\stackrel{t}{U}]{x^i}$ and choose pairwise different values $(s_1,\ldots, s_n)$ for $s$. Add the resulting equations, weighting the $i$\textsuperscript{th} equation with $(-1)^{n-1} \mu_{n-1}(\hat{s}_i) / { \prod_{i \neq j}(s_i-s_j)}$. Here $\mu_{n-1}(\hat{s}_i)$ is the elementary symmetric polynomial of degree $n-1$ in the variables $(s_1,\ldots,s_{i-1},s_{i+1},\ldots s_n)$. On the right hand side this gives the coefficient of $s^{n-1}$ of $\stackrel{s}{K}$ which is the identity operator (when considered as a (1,1)-tensor) acting on the differential of $\stackrel{t}{U}$. On the left hand side we identify the sum $\sum_{i=1}^{n}(-1)^{n-1} \mu_{n-1}(\hat{s}_i)/( \prod_{i \neq j}(s_i-s_j)) \pdiff[\stackrel{s_i}{U}]{x^i}$ with the differential of $U$ and thus arrive at $\stackrel{t}{K}\!^i_j \pdiff[U]{x^i}= \pdiff[\stackrel{t}{U}]{x^j}$.\\
For the other direction, consider two arbitrary values $s$ and $t$ and the equations
\begin{equation}\label{secondequivalencecalc}
\begin{split}
\stackrel{t}{K}\!^i_j \pdiff[U]{x^i}= \pdiff[\stackrel{t}{U}]{x^j}, \qquad 
\stackrel{s}{K}\!^i_j \pdiff[U]{x^i}= \pdiff[\stackrel{s}{U}]{x^j}
\end{split}
\end{equation}
We mulitply the first equation with $\stackrel{s}{K}\!^j_k$, and use the commutativity of $\stackrel{t}{K}$ with $\stackrel{s}{K}$ (again considered as mapping one-forms to one-forms):
\begin{equation}
\stackrel{t}{K}\!^j_k\stackrel{s}{K}\!^i_j \pdiff[U]{x^i}= \stackrel{s}{K}\!^j_k\pdiff[\stackrel{t}{U}]{x^j}
\end{equation}
Now we can use the second equation of (\ref{secondequivalencecalc}) to replace $\stackrel{s}{K}\!^i_j \pdiff[U]{x^i}$ with $\pdiff[\stackrel{s}{U}]{x^j}$ arriving back at $\stackrel{t}{K}\!^i_j \pdiff[\stackrel{s}{U}]{x^i}=\stackrel{s}{K}\!^i_j \pdiff[\stackrel{t}{U}]{x^i}$, as we desired.

\noindent To \ref{fourthequivalence}:
Whenever one applies the quantization rules \eqref{eq: Quantizationrules} to a linear polynomial $\stackrel{s}{L}$ and a polynomial $\stackrel{t}{I}$ of second degree in the momentum variables on $T^* \mathcal{M}$ and takes the commutator of the operators, then combining equations (3.8) and (3.9) from
\cite{duvalvalent2005} gives us the formula
\begin{equation}\label{ILcommutatorexpanded}
\begin{split}
[\stackrel{t}{\hat{I}}, \stackrel{s}{\hat{L}}]
= i \widehat{\{\stackrel{t}{I},  \stackrel{s}{L} \}}- \frac{i}{2} \nabla_j (\stackrel{t}{K}\!^{jk} \nabla_k (\nabla_l \stackrel{s}{V}\!^l))
\end{split}
\end{equation}
It can be obtained via explicit calculation. The second term on the right hand side of (\ref{ILcommutatorexpanded}) acts on functions by mere multiplication. It vanishes in our case because $\stackrel{t}{V}$ is a Killing vector field and thus divergence free. To show that the first term on the right hand side of (\ref{ILcommutatorexpanded}) vanishes we use that $\widehat{\{\stackrel{t}{I},  \stackrel{s}{L}\}}= - \nabla_j \circ (\{\stackrel{t}{I},  \stackrel{s}{V}\}^\sharp)^{jk} \circ \nabla_k$. So, we must show that $\{\stackrel{t}{I},  \stackrel{s}{V}\}^\sharp$ vanishes in order for $\widehat{\{\stackrel{t}{I},  \stackrel{s}{L} \}}$ to vanish. Inspection of the components of $\{\stackrel{t}{I},  \stackrel{s}{L}\}^\sharp$ reveals that they are simply the components of the Lie derivative of $\stackrel{t}{K}$ with respect to $\stackrel{t}{V}$:
\begin{equation*}
	(\{\stackrel{t}{I},  \stackrel{s}{V}\}^\sharp)^{jk} = {\stackrel{s}{V}\!^i} \partial_i {\stackrel{t}{K}\!^{jk}} - \partial_i ( {\stackrel{s}{V}\!^j} ) {\stackrel{t}{K}\!^{ik}} - \partial_i ( {\stackrel{s}{V}\!^k} ) {\stackrel{t}{K}\!^{ji}} = (\mathcal{L}_{\stackrel{s}{V}} \stackrel{t}{K})^{jk}
\end{equation*}
Applying the Leibniz rule gives 
\begin{equation*}
	\begin{split}
	(\mathcal{L}_{\stackrel{s}{V}} \stackrel{t}{K})^{jk}=&\mathcal{L}_{\stackrel{s}{V}}(\sqrt{\det (t Id - A)})\ {(tId-A)^{-1}}\!^i_l g^{lj} \\
	&+ \sqrt{\det (t Id - A)}\ (\mathcal{L}_{\stackrel{s}{V}}(tId-A)^{-1})^i_l g^{lj}\\
	&+ \sqrt{\det (t Id - A)}\ {(tId-A)^{-1}}^i_l (\mathcal{L}_{\stackrel{s}{V}}g)^{lj}
	\end{split}
\end{equation*} 
From this we see that, since the flow of $\stackrel{s}{V}$ preserves $A$ (and thus $\det A$ and functions thereof) \cite[lemma 2.2]{bolsinov2015localc} and since $\stackrel{s}{V}$ also is a Killing vector field, the term $\widehat{\{\stackrel{t}{I},  \stackrel{s}{V}\}}$ in \eqref{ILcommutatorexpanded} vanishes. Using this, a direct calculation immediately reveals that both $\{\stackrel{t}{I}+\stackrel{t}{U}, \stackrel{s}{L}\}=0 $ and $[\stackrel{t}{\hat{I}}+\stackrel{t}{\hat{U}}, \stackrel{s}{\hat{L}}]=0$
reduce to the same expression, namely $\opd \stackrel{t}{U} (\stackrel{s}{V})=0\ \forall s,t \in \mathbb{R}$.

\begin{lemma}\label{semisimplepotentiallemma}
	Let $(g,J,A)$ be c-compatible and $\stackrel{t}{K}$ be defined as in (\ref{Killingtensordefinition}). Consider a simply connected domain where the number of different eigenvalues of $A$ is constant. Let $A$ be semi-simple. Let $\stackrel{\operatorname{nc}}{E}=\{\varrho_1, \ldots, \varrho_r\}$ be the set of non-constant eigenvalues of $A$. Let $\stackrel{\operatorname{c}}{E}=\{\varrho_{r+1}, \ldots, \varrho_{r+R}\}$ be the set of constant eigenvalues and $E=\stackrel{\operatorname{nc}}{E}\cup\stackrel{\operatorname{c}}{E}$. Denote by $m(\varrho_i)$ the algebraic multiplicity of $\varrho_i$. Let  $U$ be a function such that $\stackrel{t}{K}\!^i_j \pdiff[U]{x^i}$ is exact for all values of $t$ and Let $\stackrel{t}{U}$ be such that 
	\begin{equation}\label{semisimplepotentialrequirement}
	\stackrel{t}{K}\!^i_j \pdiff[U]{x^i}= \pdiff[\stackrel{t}{U}]{x^j}
	\end{equation}
	is satisfied for all values of $t$. Then up to addition of a function of the single variable $t$ the family of functions $\stackrel{t}{U}(t,x)$
	may be written as 
	\begin{equation}\label{semisimplepotentialform1}
	\stackrel{t}{U}=\prod_{\varrho_l \in \stackrel{\operatorname{c}}{E}}^{} (t-\varrho_l)^{m(\varrho_l) /2-1}\stackrel{t}{\tilde{U}}
	\end{equation}
	where $\stackrel{t}{\tilde{U}}$ is a polynomial of degree $r-1$ in $t$. Equally $\stackrel{t}{U}$ can be written as
	\begin{equation}\label{semisimplepotentialform2}
	\stackrel{t}{U}=\sum_{i=1}^{r+R} \prod_{\varrho_l \in E\setminus\{\varrho_i\}}^{} \frac{(t-\varrho_l)^{m(\varrho_l)/2}}{(\varrho_i-\varrho_l)^{m(\varrho_l)/2}} (t-\varrho_i)^{m(\varrho_i)/2-1} f_i
	\end{equation}
	where $f_i$ are functions on $\mathcal{M}$. The functions $f_i$ may however not be chosen arbitrarily.
\end{lemma}
\textbf{Proof of lemma \ref{semisimplepotentiallemma}}: 
	Because we assumed that $A$ is semi-simple, we can factorize $\stackrel{t}{K}$ into ${\stackrel{t}{K}=\prod_{l=1}^{r+R} (t-\varrho_l) ^{m_l/2-1} \stackrel{t}{\tilde{K}}}$, with $\stackrel{t}{\tilde{K}}$ being a polynomial of degree $r+R-1$.
\begin{equation}\label{eq:55}
\prod_{\varrho_l \in \stackrel{\operatorname{c}}{E}}^{} (t-\varrho_l) ^{m(\varrho_l)/2-1} \stackrel{t}{\tilde{K}}\!^i_j \pdiff[U]{x^i}
=\pdiff[\stackrel{t}{U}]{x^i} \quad \forall t\in \mathbb{R}
\end{equation}
We used that for the non-constant eigenvalues $\varrho_1, \ldots, \varrho_r$ the multiplicities are $2$ \cite[lemma 2.2]{bolsinov2015localc}. Thus in the product on the left hand side all factors corresponding to non-constant eigenvalues are equal to $1$. We observe that upon addition of a function of the single variable $t$ to $\stackrel{t}{U}$ the equation above is still satisfied. 
This allows us to choose an arbitrary point $x_0$ and an arbitrary function $U_0 (t)$ and assume that $\stackrel{t}{U}(x_0,t)=U_0 (t)$.
Since the left hand side of \eqref{eq:55} is a polynomial in $t$ and is exact for all $t$, each of the coefficients must be exact. This allows us to integrate the terms of \eqref{eq:55} individually:
\begin{equation}\label{semisimplepotentialcomputation}
\stackrel{t}{U}(t,x)=U_0(t)+\int_{x_0}^{x} \stackrel{t}{K}(\opd U)=U_0(t)+\prod_{\varrho_l \in \stackrel{\operatorname{c}}{E}} (t-\varrho_l) ^{m(\varrho_l)/2-1}\sum_{i=0}^{r} t^i \int_{x_0}^{x}\tilde{K}_{(i)}(\opd U)
\end{equation}
The last step of this calculation makes use of the fact that if $\varrho_l$ is of multiplicity $m_l \geq 4$ then $\varrho_l$ is constant \cite[lemma 2.2]{bolsinov2015localc}.
The integral is meant to be taken along any path connecting $x_0$ and $x$ and the $\tilde{K}_{(i)}$ is the coefficient of $t^i$ in $\stackrel{t}{\tilde{K}}$. Again $\stackrel{t}{\tilde{K}}$ and $\tilde{K}_{(i)}$ are considered as $(1,1)$ tensors mapping 1-forms to 1-forms. So for any value of $t$ the value of $\stackrel{t}{U}$ at $x$ is uniquely defined by its value at $x_0$ and the function $U$. Formula (\ref{semisimplepotentialcomputation}) proves the claim that $\stackrel{t}{U}$ can be written in the form (\ref{semisimplepotentialform1}) where on the right hand side $U_0(t)$ takes the role of the possible addition of a function of $t$ alone. Evidently, we have $\stackrel{t}{\tilde{U}}(x)=\sum_{i=0}^{r+R} t^i \int_{x_0}^{x}\tilde{K}_{(i)}(\opd U)$. Since $\stackrel{t}{\tilde{U}}$ is a polynomial of degree $r+R-1$ it is uniquely defined by its values at the $r+R$ different eigenvalues of $A$. Via the Lagrange interpolation formula we have 

\begin{equation}
\stackrel{t}{\tilde{U}}=\sum_{i=1}^{r+R} \prod_{\varrho_l \in E\setminus\{\varrho_i\}}^{} \frac{(t-\varrho_l)}{(\varrho_i-\varrho_l)} \tilde{f}_i
\end{equation}
for some funcions $\tilde{f}_i$. 
Introducing $f_i=\prod_{\varrho_l \in E\setminus\{\varrho_i\}} (\varrho_i-\varrho_l)^{m(\varrho_l)/2 -1} \tilde{f}_i$ the potential $U$ can be written as
\begin{equation}
\stackrel{t}{U}=\sum_{i=1}^{r+R} \prod_{\varrho_l \in E\setminus\{\varrho_i\}}^{} \frac{(t-\varrho_l)^{m(\varrho_l)/2}}{(\varrho_i-\varrho_l)^{m(\varrho_l)/2}} (t-\varrho_i)^{m(\varrho_i)/2-1} f_i
\end{equation}
concluding the proof of lemma \ref{semisimplepotentiallemma}.

\begin{lemma}\label{functiondifferentials are eigenforms}
	Let $(g,J,A)$ be c-compatible and $\stackrel{t}{K}$ as in (\ref{Killingtensordefinition}). Let $A$ be semi-simple. Let $\stackrel{\operatorname{nc}}{E}=\{\varrho_1, \ldots, \varrho_r\}$ be the set of non-constant eigenvalues of $A$. Let $\stackrel{\operatorname{c}}{E}=\{\varrho_{r+1}, \ldots, \varrho_{r+R}\}$ be the set of constant eigenvalues and $E=\stackrel{\operatorname{nc}}{E}\cup\stackrel{\operatorname{c}}{E}$. The multiplicity of $\varrho_l$ is denoted by $m(\varrho_l)$. Let
	\begin{equation}
	\stackrel{t}{U}=\sum_{i=1}^{r+R} \prod_{\varrho_l \in E\setminus\{\varrho_i\}}^{} \frac{(t-\varrho_l)^{m(\varrho_l)/2}}{(\varrho_i-\varrho_l)^{m(\varrho_l)/2}} (t-\varrho_i)^{m(\varrho_i)/2-1} f_i
	\end{equation}
	and let 
	\begin{equation}\label{eigenformrequirement}
	\stackrel{t}{K}\!^i_j \pdiff[U]{x^i}= \pdiff[\stackrel{t}{U}]{x^j}
	\end{equation}
	be satisfied for all values of $t$. Then for all values of $i$ the relation $\opd f_i \circ A = \varrho_i \opd f_i$ must be satisfied. In other words:  The differentials of the functions $f_i$ are eigenvectors of $A$ with eigenvalues $\varrho_i$, where $A$ is considered as to map one-forms to one-forms.
\end{lemma}
\textbf{Proof of lemma \ref{functiondifferentials are eigenforms}}: We consider $\stackrel{t}{K}$ as a $(1,1)$ tensor field.
Using our assumption that $A$ is semi-simple we rewrite equation \eqref{eigenformrequirement} in terms of the quantities $\stackrel{t}{\tilde{U}}$, $\tilde{f_i}$ and $\stackrel{t}{\tilde{K}}$ defined by
\begin{align*}
	\stackrel{t}{\tilde{U}}&=\sum_{i=1}^{r+R} \prod_{\varrho_l \in E\setminus\{\varrho_i\}}^{} \frac{(t-\varrho_l)}{(\varrho_i-\varrho_l)} \tilde{f}_i\\
	f_i &=\prod_{\varrho_l \in E\setminus\{\varrho_i\}}^{} (\varrho_i-\varrho_l)^{m(\varrho_l)/2 -1} \tilde{f}_i\\
	\stackrel{t}{K}&=\prod_{\varrho_l \in E\setminus\{\varrho_i\}} (t-\varrho_l) ^{m(\varrho_l)/2-1} \stackrel{t}{\tilde{K}}
\end{align*}
Then we use the fact that eigenvalues of multiplicity $m_l \geq 4$ are constant \cite[lemma 2.2]{bolsinov2015localc} and equation \eqref{eigenformrequirement} transforms into
\begin{equation}\label{eq:61}
\stackrel{t}{\tilde{K}} (\opd U) = \opd\stackrel{t}{\tilde{U}}
\end{equation}
by dividing out the common factors.\\
The right hand side can be rewritten: consider $\stackrel{t}{\tilde{U}}$ where, rather than choosing a constant value for the parameter $t$ we fill in the $l$\textsuperscript{th} eigenvalue of $A$. Then we have $\stackrel{\varrho_l}{\tilde{U}}= \tilde{f}_l$. Taking the differential and rearranging the terms gives
\begin{equation}\label{eq:62}
\opd \tilde{f}_l - \left.\pdiff[\stackrel{t}{\tilde{U}}]{t}\right|_{t=\varrho_l} \opd \varrho_l =\opd \stackrel{\varrho_l}{\tilde{U}} - \left.\pdiff[\stackrel{t}{\tilde{U}}]{t}\right|_{t=\varrho_l} \opd \varrho_l = \opd \stackrel{t}{\tilde{U}}|_{t=\varrho_l}
\end{equation}
We evaluate \eqref{eq:61} at $t=\varrho_l$ and plug in \eqref{eq:62}:
\begin{equation}\label{eq:87}
\stackrel{\varrho_l}{\tilde{K}} (\opd U)= \opd \tilde{f}_l - \left.\pdiff[\stackrel{t}{\tilde{U}}]{t}\right|_{t=\varrho_l} \opd \varrho_l
\end{equation}
Because we assumed $A$ to be semi-simple we can decompose $\opd U$ into  one-forms~$\upsilon_l$  such that $\upsilon_l \circ A = \varrho_l \upsilon_l$. 
From the definition of $\stackrel{t}{\tilde{K}}$ we have that $
\stackrel{t}{\tilde{K}}\!(\upsilon_l)= \prod_{\varrho_m \in E\setminus\{\varrho_l\}}^{}(t- \varrho_m) \upsilon_l$, again because we assumed $A$ to be semi-simple. Evaluating at $t=\varrho_k$ yields
\begin{equation}\label{eq:88}
\begin{split}
\stackrel{\varrho_k}{\tilde{K}}(\upsilon_l) = \left(\prod_{\varrho_m \in E\setminus\{\varrho_l\}} (\varrho_k -\varrho_m) \right) \upsilon_l
\end{split}
\end{equation}
In particular this means that if $k\neq l$ then $\stackrel{\varrho_k}{\tilde{K}}(\opd \varrho_l)$ is zero.
Plugging this into \ref{eq:87} we get that on the left hand side only $\stackrel{\varrho_l}{\tilde{K}}(\opd U)=\stackrel{\varrho_l}{\tilde{K}}(\upsilon_l)$, which we express via \eqref{eq:88}:
\begin{equation}
\left(\prod_{\varrho_m \in E\setminus\{\varrho_l\}} (\varrho_l -\varrho_m) \right) \upsilon_l= \opd \tilde{f}_l - \left.\pdiff[\stackrel{t}{\tilde{U}}]{t}\right|_{t=\varrho_l} \opd \varrho_l
\end{equation}
Since $\opd \varrho_l \circ A = \varrho_l \opd \varrho_l$ \cite[lemma 2.2]{bolsinov2015localc} and $\upsilon_l \circ A = \varrho_l \upsilon_l$ we have that $\opd \tilde{f}_l \circ A = \varrho_l \opd \tilde{f}_l$. The way in which $\tilde{f}_l$ was constructed  then implies $\opd f_l \circ A = \varrho_l \opd f_l$, concluding the proof of lemma \ref{functiondifferentials are eigenforms}.
\begin{lemma}\label{lemma:24}
	Let $(g,J,A)$ be c-compatible, $\stackrel{t}{K}$ as in \eqref{Killingtensordefinition} and $\stackrel{t}{V}$ as in \eqref{Killingvectorfielddefinition}.  $\stackrel{\operatorname{nc}}{E}=\nobreak \{\varrho_1, \ldots, \varrho_r\}$ denotes the set of non-constant eigenvalues of $A$.  $\stackrel{\operatorname{c}}{E}=\{\varrho_{r+1}, \ldots, \varrho_{r+R}\}$ denotes the set of constant eigenvalues and $E=\stackrel{\operatorname{nc}}{E}\cup\stackrel{\operatorname{c}}{E}$. The multiplicity of $\varrho_l$ is denoted by $m(\varrho_l)$. Let 
	\begin{equation}\label{eq:90}
	\stackrel{t}{U}=\sum_{i=1}^{r+R} \prod_{\varrho_l \in E\setminus\{\varrho_i\}}^{} \frac{(t-\varrho_l)^{m(\varrho_l)/2}}{(\varrho_i-\varrho_l)^{m(\varrho_l)/2}} (t-\varrho_i)^{m(\varrho_i)/2-1} f_i
	\end{equation}
	Let furthermore $\opd f_l \circ A = \varrho_l \opd f_l$ for all values of $l=1 \ldots r+R$.
	Let $\opd \stackrel{t}{U} (\stackrel{s}{V})=0$ be satisfied for all values of $s,t \in \mathbb{R}$.
	Then $\opd \stackrel{t}{U} (\stackrel{s}{V})=0$ is satisfied for all values of $s,t \in \mathbb{R}$ if and only if for each eigenvalue $\varrho_k$ of $A$ the differential $\opd f_k$ is proportional to $\opd \varrho_k$ at all points where $\opd \varrho_k \neq 0$.
\end{lemma}
\begin{corollary}
	If $\varrho_l$ is a non-constant real eigenvalue and $\opd\varrho_l\neq 0$ in the neighbourhood of a given point then locally $f_l$ can be expressed as a smooth function of $\varrho_l$. Likewise if $\varrho_l$ is a non-constant complex eigenvalue and $\opd\varrho_l\neq 0$ in the neighbourhood of a given point then locally $f_l$ can be expressed as a holomorphic function of $\varrho_l$.
\end{corollary}
\textbf{Proof of lemma \ref{lemma:24}}:
The condition $\opd \stackrel{t}{U} (\stackrel{s}{V})=0\ \forall\, s,t\in \mathbb{R}$ is equivalent to \\$\opd\! \stackrel{t}{U}\!(\operatorname{span}\{ J \grad \varrho_i, i\!=\!1\ldots r\})\!=\!0$, because $\operatorname{span}\{\stackrel{t}{V},t\in \mathbb{R}\}\!=\! \operatorname{span}\{J \grad \varrho_i, {i\!=\!1\ldots r}\}$.
From \eqref{eq:90} we see that $\opd \stackrel{t}{U}$ involves (with some coefficients) only the differentials of the eigenvalues of $A$ and the differentials of the functions $f_i$. Thus $\opd \stackrel{t}{U}(J \grad \varrho_i)$ is a linear combination of $\opd \varrho_j$ and $\opd f_j$ applied to $J \grad \varrho_i$.
But $\opd \varrho_j (J \grad \varrho_i)$ is zero for all values of $i,j$: if $i=j$, then $\opd \varrho_j (J \grad \varrho_i)=0$ due to the fact that $J$ is antisymmetric with respect to $g$ and if $i\neq j$ then $\opd \varrho_j (J \grad \varrho_i)=0$, because $A$ is $g$-self-adjoint and $\grad \varrho_i$ and $\grad \varrho_j$ are eigenvectors of $A$ with different eigenvalues. 
Thus $\opd \stackrel{t}{U}(J \grad \varrho_i)$ is a linear combination of $\{\opd f_j (J \grad \varrho_i), j=1\ldots r+R\}$.
But because we assumed that for all $l$ $\opd f_l \circ A = \varrho_l \opd f_l$ and because $A \grad \varrho_i = \varrho_i \grad \varrho_i$ \cite[lemma2.2]{bolsinov2015localc} and $A$ commutes with $J$ and is $g$-self-adjoint we get that $\opd \stackrel{t}{U}(J\grad \varrho_i)$ is some coefficient times $\opd f_i(J \grad\varrho_i)$.\\
From \eqref{eq:90} we see that this coefficient is 
\begin{equation*}
	\prod_{\varrho_l \in E\setminus\{\varrho_i\}}^{} \frac{(t-\varrho_l)^{m(\varrho_l)/2}}{(\varrho_i-\varrho_l)^{m(\varrho_l)/2}} (t-\varrho_i)^{m(\varrho_i)/2-1}
\end{equation*} 
But for a given value of $t$ this can only vanish at points on $\mathcal{M}$ where $t$ is equal to an eigenvalue of $A$. So at each point on the manifold we can choose a value for $t$ such that this expression is non-zero and thus $\opd f_i (J\grad \varrho_i)$ must vanish for all values of $i$. 
If we consider a value $i$ such that $\varrho_i$ is constant then $\opd f_i (J\grad \varrho_i)=0$ is trivially satisfied. If $\varrho_i$ is non-constant then its multiplicity is $2$ \cite[lemma 2.2]{bolsinov2015localc} and at all points where $\opd \varrho_i \neq 0$ the set $\{\grad \varrho_i, J\grad \varrho_i\}$ is an orthogonal basis of the $\varrho_i$-eigenspace of $A$. It follows that at such points $\opd f_i$ may be written as a linear combination of $\opd \varrho_i$ and $\opd \varrho_i \circ J$. Plugging this decomposition into $\opd f_i (J\grad \varrho_i)=0$ we conclude that $\opd f_i$ is proportional to $\opd \varrho_i$ at all points where $\opd \varrho_i\neq 0$ and lemma \ref{lemma:24} is proven.

\begin{lemma}\label{lemma:25}
	Let $(g,J,A)$ be c-compatible.  $\stackrel{\operatorname{nc}}{E}=\nobreak \{\varrho_1, \ldots, \varrho_r\}$ denotes the set of non-constant eigenvalues of $A$.  $\stackrel{\operatorname{c}}{E}=\{\varrho_{r+1}, \ldots, \varrho_{r+R}\}$ denotes the set of constant eigenvalues and $E=\stackrel{\operatorname{nc}}{E}\cup\stackrel{\operatorname{c}}{E}$. The multiplicity of $\varrho_l$ is denoted by $m(\varrho_l)$. Let
	\begin{equation*}
	\stackrel{t}{K}\!^{ij} \stackrel{\operatorname{def}}{=}\sqrt{\det (t Id - A)}\ {(tId-A)^{-1}}^i_l g^{lj}, \qquad \stackrel{t}{V}\!^j = J^j_k g^{ki} \nabla_i \sqrt{\det(t Id -A)}
	\end{equation*}
	and 
	\begin{equation}
	\stackrel{t}{U}=\sum_{i=1}^{r+R} \prod_{\varrho_l \in E\setminus\{\varrho_i\}}^{} \frac{(t-\varrho_l)^{m(\varrho_l)/2}}{(\varrho_i-\varrho_l)^{m(\varrho_l)/2}} (t-\varrho_i)^{m(\varrho_i)/2-1} f_i
	\end{equation}
	with $\opd f_l \circ A = \varrho_l \opd f_l$ for all $l=1\ldots r$ and $\opd f_l$ proportional to $\opd \varrho_l$ for all $l$ for which $\varrho_l$ is non-constant.\\
	Then
	\begin{equation}
	\stackrel{t}{K}\!^i_j \pdiff[\stackrel{s}{U}]{x^i} = \stackrel{s}{K}\!^i_j \pdiff[\stackrel{t}{U}]{x^i}\ \forall s,t \in \mathbb{R}\quad and \quad \opd \stackrel{t}{U} (\stackrel{s}{V})=0\ \forall s,t \in \mathbb{R}
	\end{equation}
\end{lemma}
\textbf{Proof of lemma \ref{lemma:25}}: 
$\opd \stackrel{t}{U} (\stackrel{s}{V})=0\ \forall s,t \in \mathbb{R}$ is fulfilled because $\opd f_i$ and $\opd \varrho_i$ evaluate to zero when applied to $J \grad \varrho_j$ for all $i,j=1,\ldots, r+R$.

To see $\stackrel{t}{K}\!^i_j \pdiff[\stackrel{s}{U}]{x^i} = \stackrel{s}{K}\!^i_j \pdiff[\stackrel{t}{U}]{x^i}$, we compute $\opd \stackrel{t}{U}$ using the fact that non-constant eigenvalues of $A$ are of multiplicity 2 \cite[lemma 2.2]{bolsinov2015localc}:
\begin{equation}\label{eq:99}
\begin{split}
\opd \stackrel{t}{U}=& \sum_{i=1}^{r+R} \prod_{\varrho_l \in E\setminus\{\varrho_i\}}\!\! \left( \frac{t-\varrho_l}{\varrho_i-\varrho_l}\right)^{m(\varrho_l)/2} \left(t-\varrho_i\right)\!^{m(\varrho_i)/2-1} \left[ \opd f_i - \sum_{\mathclap{\varrho_p\in E\setminus\{\varrho_i\}}}^{}\ \  \frac{m(\varrho_p)/2}{\varrho_i -\varrho_p} f_i \opd \varrho_i\right]\\
& - \sum_{i =1}^{r+R} \prod_{\varrho_l \in E\setminus\{\varrho_k\}}^{} \frac{t- \varrho_l}{(\varrho_i - \varrho_l)^{m(\varrho_l)/2}} \sum_{\varrho_p \in \stackrel{\operatorname{nc}}{E}\setminus\{\varrho_i\}}^{} \prod_{\varrho_l \in E\setminus\{\varrho_p\}}^{}(t- \varrho_l)^{m(\varrho_l)/2}  f_i \frac{\opd \varrho_p}{\varrho_p- \varrho_i}
\end{split}
\end{equation}
Considering $\stackrel{s}{K}$ as a $(1,1)$-tensor mapping one-forms to one-forms and using that for all $i=1 \ldots r$: $\opd \varrho_i \circ A = \varrho_i \opd \varrho_i$ and $\opd f_i \circ A = \varrho_i \opd f_i$ we have 
\begin{equation}\label{eq:100}
\begin{split}
	\stackrel{s}{K}(\opd \varrho_i) &=\bigg(\prod_{\varrho_l \in E\setminus\{\varrho_i\}}^{} (s-\varrho_l)^{m_l/2}\bigg) (s- \varrho_i)^{m_i/2 -1} \opd \varrho_i \\
	\stackrel{s}{K}(\opd f_i) &=\bigg(\prod_{\varrho_l \in E\setminus\{\varrho_i\}}^{} (s-\varrho_l)^{m_l/2}\bigg) (s- \varrho_i)^{m_i/2 -1} \opd f_i
\end{split}
\end{equation}
Again we consider $\stackrel{s}{K}$ as a $(1,1)$-tensor acting on the differential of $\stackrel{t}{U}$. By combining (\ref{eq:99}) and (\ref{eq:100}) and using that the non-constant eigenvalues have multiplicity $2$ we get:
\begin{equation}
\begin{split}
\stackrel{s}{K}(\opd \stackrel{t}{U})=& \sum_{i=1}^{r+R} \prod_{\varrho_l \in E\setminus\{\varrho_i\}} \left(\frac{(s-\varrho_l)(t-\varrho_l)}{\varrho_i-\varrho_l}\right)^{m_l/2} \left(\left(t-\varrho_i\right)\left(s-\varrho_i\right)\right)^{m_i/2-1}\\
&\phantom{+}\times \left[ \opd f_i - \sum_{p\neq i}^{} \frac{m_p/2}{\varrho_i -\varrho_p} f_i \opd \varrho_i\right]\\
& - \sum_{i =1}^{r+R} \prod_{\varrho_l \in E\setminus\{\varrho_k\}}^{} \frac{1}{(\varrho_i - \varrho_l)^{m_l/2}}\\
&\phantom{+}\times \hspace{-1em}\sum_{\varrho_p \in \stackrel{\operatorname{nc}}{E}\setminus\{\varrho_i\}}^{} \prod_{\varrho_l \in E\setminus\{\varrho_p\}}^{}\left((s- \varrho_l)(t- \varrho_l)\right)^{m_l/2}  f_i \frac{\opd \varrho_p}{\varrho_p- \varrho_i}
\end{split}
\end{equation}
The right hand side is apparently symmetric when exchanging $s$ and $t$ and as a consequence $\stackrel{t}{K}\!^i_j \pdiff[\stackrel{s}{U}]{x^i}=\stackrel{s}{K}\!^i_j \pdiff[\stackrel{t}{U}]{x^i}$  is fulfilled, concluding the proof of lemma \ref{lemma:25}.

\textbf{Proof of theorem \ref{KommutierendePotentiale}}: The theorem results as a combination of the proofs of lemmata \ref{lemma:4problems} and \ref{lemma:25}.

\textbf{Proof of theorem \ref{PotentialForm}}: combining the proofs of lemmata \ref{lemma:4problems}, \ref{semisimplepotentiallemma}, \ref{functiondifferentials are eigenforms}, \ref{lemma:24} and $\ref{lemma:25}$.

\subsection{Proof of theorem \ref{eigenfunctionstheorem} | Common Eigenfunctions}\label{commoneigenfunctions}
We first show that $\psi$ is an eigenfunction of $\stackrel{s}{\hat{L}}$ for all values of $s$, if and only if it is an eigenfunction of $\partial_{t_i}$ for all values of $i$. 
A direct computation from \eqref{eq:example general 1}, \eqref{eq:example general 2}, \eqref{eq:example general 3} shows:
\begin{equation}
	\stackrel{s}{\hat{L}}\psi = \prod_{\varrho_p \in \stackrel{\operatorname{c}}{E}}^{}(s-\varrho_p)^{m(\varrho_p)}\sum_{q=1}^{r} s^{r-q} \partial_{t_q} \psi
\end{equation}
Now suppose that $\psi$ is an eigenfunction of $\stackrel{s}{\hat{L}}$ for all values of $s$ and denote the eigenvalue by $\stackrel{s}{\omega}$. Then $\stackrel{s}{\omega}$ must be a polynomial in $s$ of degree $r+R-1$, because $\stackrel{s}{\hat{L}}$ is a polynomial in $s$ of degree $r+R-1$. From the equation above we see that $\stackrel{s}{\omega}$ must have a zero of order $m(\varrho_p)$ at $s=\varrho_p$ for all constant eigenvalues $\varrho_p$. Thus we can write $\stackrel{s}{\omega}=\prod_{\varrho_p \in \stackrel{\operatorname{c}}{E}}^{}(s-\varrho_p)^{m(\varrho_p)}\sum_{q=1}^{r} s^{r-q} i \omega_q$. But because polynomials are equal if and only if all their coefficients are equal we get that $\stackrel{s}{\hat{L}} \psi = \stackrel{s}{\omega} \forall s \in \mathbb{R}$ if and only if $i \partial_{t_q} \psi = \omega_q \psi$ for $q=1, \ldots, r$.

\noindent To obtain the other separated equations we work with the family of second order differential operators $\stackrel{s}{K}$. Since the metric is not given in terms of the coordinate basis and the one-forms $\alpha$ and $\vartheta$ are not unique it poses an obstruction to using the standard formula for the Laplacian and $\stackrel{t}{\hat{K}}$. The workaround is quick and simple though.

\noindent Let $\{X_i, i=1\ldots n\}$ be a set of $n$ linearly independent differentiable vector fields on a manifold $\mathcal{M}^n$ and denote by $\{\beta^i, i=1\ldots n\}$ its dual basis, i.e.\ $\beta^i (X_j) = \delta^i_j$. We shall denote by $T$ the matrix relating the coordinate Vector fields $\partial_i$ and the vector fields $X_i$: $T^i_j X_i = \partial_j$
Then for an arbitrary symmetric $(2,0)$-tensor the following formula is easily obtained via the product rule for partial derivatives:
\begin{equation}\label{eq:pseudolaplacianinarbitraryframe}
	\begin{split}
	\frac{1}{\sqrt{\left|\det g \right|}} \partial_i \sqrt{\left|\det g \right|} \stackrel{s}{K}\!^{ij} \partial_j =& \frac{\det T}{\sqrt{\left| \det g \right|}} X_s \frac{\sqrt{\left| \det g \right|}}{\det T} T^s_i \stackrel{s}{K}\!^{ij} T^k_j X_j \\
	&- X_s(T^s_i) \stackrel{s}{K}\!^{ij} T^k_j X_k + \frac{X_s(\det T)}{\det T} T^s_i \stackrel{s}{K}\!^{ij} T^k_j X_j
	\end{split}
\end{equation}
where on the right hand side the $X$'s are to be interpreted as the directional derivative in the sense $X_s = (T^{-1})^i_s \partial_i$ and in the last two terms $X_s(\cdot)$ is meant as to only act on the expression in the parentheses.

\noindent The quantity $\frac{\det g}{\det T^2}$ is simply the determinant of the matrix with $(i,j)$\textsuperscript{th} component $g(X_i,X_j)$ and $T^s_i \stackrel{s}{K}\!^{ij} T^k_j$ are the components of $\stackrel{s}{K}\!$ in the basis $\{X_i, i=1\ldots n\}$. If $\{X_i, i=1\ldots n\}$ are the coordinate vector fields belonging to some coordinate system then the last two terms cancel out and one arrives at the well known fact that the left hand side expression is independent of the choice of coordinates. 

\noindent In our case we choose $(X_i) = (\partial_{\chi_i}, \partial_{t_i}, \partial_{\stackrel{\gamma}{y}_i} - \sum_{p=1}^{r} \stackrel{\gamma}{\alpha}_{pi} \partial_{t_p})$. The dual basis consists of the one-forms $(\opd \chi_i, \opd t_i + \alpha_i , \opd \stackrel{\gamma}{y}_i)$. Taking $i$ as the column index and $j$ as the row index we have the components $T^i_j$ given by
\begin{equation*}
	\left(
	\begin{array}{c|c|c}
	Id 		& 0			& 0\\
			\hline
	0		& Id 		& 0	\\
			\hline
	0		& *			& Id
	
	\end{array}
	\right)
\end{equation*}
Where the $*$-block contains the components of $\alpha_l$ as the $l$\textsuperscript{th} column and the $Id$-Blocks are of dimensions equal to the number of $\chi$-, $t$- and $y$-coordinates. 

\noindent From this we can conclude that $X_s (T^s_i)= (T^{-1})^j_s \partial_j T^s_i=0$ because the one-forms $\alpha$ do not depend on  the $t$-variables. 
Furthermore $\det T=1$ and thus $X_s(\det T)=0$ for all values of $s$.

\noindent For our specific case \eqref{eq:pseudolaplacianinarbitraryframe} simplifies to
\begin{equation}
	\begin{split}
	\frac{1}{\sqrt{\left|\det g \right|}} \partial_i \sqrt{\left|\det g \right|} \stackrel{s}{K}\!^{ij} \partial_j =& \frac{\det T}{\sqrt{\left| \det g \right|}} X_s \frac{\sqrt{\left| \det g \right|}}{\det T} T^s_i \stackrel{s}{K}\!^{ij} T^k_j X_j \\
	\end{split}
\end{equation}

\noindent Here $\frac{\sqrt{|\det g|}}{\det T}$ is the determinant of the matrix of $g$ in the basis $(\opd \chi_i, \opd t_i, \opd \stackrel{\gamma}{y_i})$ and $T^s_i K^{ij} T^k_j$ is the matrix of $\stackrel{s}{K}\!$ in the basis $(X_i)$. These quantities can be obtained from the formulae \eqref{eq:example general 1}, \eqref{eq:example general 2} and \eqref{eq:example general 3}. 
We can then express the Vector fields $(X_i)$ in terms of the coordinate basis to get the following result (We use that $A_\gamma = \varrho_{\gamma} Id$, because we assumed that A is semi-simple and that all constant eigenvalues are real.):
\begin{align}\label{eq:pseudolaplacian in explicit coordinates}
\begin{aligned}\raisetag{18em}
	&\nabla_i \stackrel{s}{K}\!^{ij} \nabla_j = \frac{1}{\sqrt{\left|\det g \right|}} \partial_i \sqrt{\left|\det g \right|} \stackrel{s}{K}\!^{ij} \partial_j =\\
	& \hspace{0cm}\sum_{\varrho_k \in \stackrel{\operatorname{nc}}{E}}^{} \frac{\prod_{\varrho_l \in \stackrel{\text{}}{E}\setminus\{\varrho_k\}} (s-\varrho_l)^{m(\varrho_l)/2}}{\varepsilon_k \Delta_k \prod_{\varrho_\gamma \in \stackrel{\operatorname{c}}{E}}^{} (\varrho_\gamma - \varrho_k)^{m(\varrho_\gamma)/2} \varrho_k'} \partial_{\chi_k} \varrho_k' \prod_{\varrho_\gamma \in \stackrel{\operatorname{c}}{E}} (\varrho_\gamma - \varrho_k)^{m(\varrho_\gamma)/2} \partial_{\chi_k}\\
	& \hspace{0cm}+ \sum_{i,j=1}^{r} \sum_{\varrho_k \in \stackrel{\operatorname{nc}}{E}}^{} \frac{\varepsilon_k (- \varrho_k)^{2r-i-j}}{\Delta_k (\varrho_k ')^2} \prod_{\varrho_l \in E\setminus\{\varrho_k\}} (s- \varrho_l)^{m(\varrho_l)/2} \partial_{t_i} \partial_{t_j}\\
	& \hspace{0cm}+\!\! \sum_{\gamma:\varrho_\gamma \in \stackrel{\operatorname{c}}{E}}^{} \frac{1}{\sqrt{|\det g_\gamma |}} \partial_{\stackrel{\gamma}{y}_i} \sqrt{|\det g_\gamma|} \bigg(\prod_{\varrho_l\in E}(s-\varrho_l)^{m(\varrho_l)/2} (s Id -A_\gamma)^{-1} \\
	&\!\! \hspace{5.6cm}\prod_{\varrho_k \in \stackrel{\operatorname{nc}}{E}} (A_\gamma - \varrho_k Id)^{-1} g_\gamma^{-1}\bigg)^{ij} \partial_{\stackrel{\gamma}{y}_j}\\
	& \hspace{0cm}-\!\! \sum_{\gamma:\varrho_\gamma \in \stackrel{\operatorname{c}}{E}}^{} \frac{1}{\sqrt{|\det g_\gamma |}} \partial_{\stackrel{\gamma}{y}_i} \sqrt{|\det g_\gamma|} \bigg(\prod_{\varrho_l\in E}(s-\varrho_l)^{m(\varrho_l)/2} (s Id -A_\gamma)^{-1}\\
	&\!\! \hspace{5.6cm}\prod_{\varrho_k \in \stackrel{\operatorname{nc}}{E}} (A_\gamma - \varrho_k Id)^{-1} g_\gamma^{-1}\bigg)^{ij} \sum_{q=1}^{r}\stackrel{\gamma}{\alpha}_{qj} \partial_{t_q}\\
	& \hspace{0cm}-\!\! \sum_{\gamma:\varrho_\gamma \in \stackrel{\operatorname{c}}{E}}^{} \! \bigg(\prod_{\varrho_l\in E}(s-\varrho_l)^{m(\varrho_l)/2} (s Id -A_\gamma)^{-1} \prod_{\varrho_k \in \stackrel{\operatorname{nc}}{E}} (A_\gamma - \varrho_k Id)^{-1} g_\gamma^{-1}\bigg)^{ij} \sum_{q=1}^{r}\stackrel{\gamma}{\alpha}_{qi} \partial_{t_q} \partial_{\stackrel{\gamma}{y}_j} \\
	& \hspace{0cm}-\!\! \sum_{\gamma:\varrho_\gamma \in \stackrel{\operatorname{c}}{E}}^{} \! \bigg(\prod_{\varrho_l\in E}(s-\varrho_l)^{m(\varrho_l)/2}(s Id -A_\gamma)^{-1} \prod_{\varrho_k \in \stackrel{\operatorname{nc}}{E}} (A_\gamma - \varrho_k Id)^{-1} g_\gamma^{-1}\bigg)^{ij} \! \sum_{p,q=1}^{r} \!\!\stackrel{\gamma}{\alpha}_{qi}\stackrel{\gamma}{\alpha}_{pj} \partial_{t_q} \partial_{t_p} 
\end{aligned}
\end{align}

If now we suppose that $\psi$ simultaneously satisfies the family of eigenvalue equations
\begin{equation}\label{eq:Eigenwertgleichung fuer s}
	- \nabla_i \stackrel{s}{K}\!^{ij} \nabla_j \psi + \stackrel{s}{U} \psi = \stackrel{s}{\lambda} \psi
\end{equation}
then the left hand side is a polynomial in $s$ and for the equations to hold $\stackrel{s}{\lambda}$ must be a polynomial of degree $n-1$ in $s$ as well and we shall write $\stackrel{s}{\lambda}=\sum_{l=0}^{n-1} \lambda_l s^l$. Because of our assumption that all constant eigenvalues are real and that $A$ is semi-simple we have $A_\gamma= c_\gamma Id= \varrho_{r+\gamma} Id$ for $\gamma=1 \ldots R$. Thus for all values of $\gamma$ where $m(\varrho_{\gamma})>2$ the left hand side of $-\nabla_j \stackrel{s}{K}\!^{jk} \nabla_k \psi + \stackrel{s}{U} \psi = \stackrel{s}{\lambda} \psi$ has a zero of multiplicity $m(\varrho_{\gamma})/2-1$ at $s=\varrho_{\gamma}$.
We can thus write $\stackrel{s}{\lambda}=\prod_{\varrho_{\gamma}\in \stackrel{\operatorname{c}}{E}}^{} (s-\varrho_{\gamma})^{m(\varrho_{\gamma})/2-1}\stackrel{s}{\tilde{\lambda}}$ where $\stackrel{s}{\tilde{\lambda}}$ is a polynomial of degree ${r+R-1}$ in $s$. $\stackrel{s}{\tilde{\lambda}}$ can be written as $\stackrel{s}{\tilde{\lambda}}=\sum_{j=0}^{r+R-1} s^j \tilde{\lambda}_j$.
Because we assumed that the eigenvalue equation \eqref{eq:Eigenwertgleichung fuer s} is satisfied for all values of $s$ we can split \eqref{eq:Eigenwertgleichung fuer s} into $n$ equations, one for each coefficient.
We choose a non-constant eigenvalue $\varrho_k \in \stackrel{\operatorname{nc}}{E}$ and multiply the equation for the coefficient of $s^l$ with $\varrho_k^l$. We do this for all $n$ equations and add them up. 
If we use $\nabla_i \stackrel{s}{K}\!^{ij} \nabla_j=\frac{1}{\sqrt{\left|\det g \right|}} \partial_i \sqrt{\left|\det g \right|} \stackrel{s}{K}\!^{ij} \partial_j$ in \eqref{eq:Eigenwertgleichung fuer s} then the result of these steps can be is simply that we replace $\nabla_i \stackrel{s}{K}\!^{ij} \nabla_j \psi$ via \eqref{eq:pseudolaplacian in explicit coordinates} and substitute $s$ with $\varrho_k$.
We end up with
\begin{multline}\label{eq:separatedode}
	\frac{-1}{\varepsilon_k \varrho_k '}\partial_{\chi_k} \varrho_k' \prod_{\varrho_{\gamma} \in \stackrel{\operatorname{c}}{E}} (\varrho_{k} - \varrho_{\gamma})^{m(\varrho_{\gamma})/2} \partial_{\chi_k} \psi \\
	- \sum_{i,j=1}^{r} \frac{\varepsilon_k (-\varrho_k)^{2r-i-j}}{(\varrho_k')^2}\prod_{\varrho_{\gamma}\in \stackrel{\operatorname{c}}{E}}^{}(\varrho_k - \varrho_{\gamma})^{m(\varrho_{\gamma})/2} \partial_{t_i} \partial_{t_j}
	+ f_r \psi = \sum_{i=0}^{n-1} \lambda_i \varrho_k^i
\end{multline}
Because $\partial_{t_q}\psi=i \omega_q \psi$, $\psi$ satisfies the ordinary differential equations that have been claimed.

\noindent To obtain the separated partial differential equations that we have claimed we divide $-\nabla_j {\stackrel{s}{K}\!^{jk}} \nabla_k \psi + \stackrel{s}{U} \psi = \stackrel{s}{\lambda} \psi$ by $\prod_{\varrho_{\gamma}\in \stackrel{\operatorname{c}}{E}}^{} (s-\varrho_{\gamma})^{m(\varrho_{\gamma})/2-1}$. Then we choose a constant eigenvalue $\varrho_{\gamma}$ and  evaluate the result at $s=\varrho_{\gamma}\in \stackrel{\operatorname{c}}{E}$ in the same way as before:
\begin{equation}\label{eq:separatedpde}
	\begin{split}
	\sum_{j=0}^{r+R-1} \tilde{\lambda}_j \varrho_{\gamma}^j
	=&-\,\prod_{\mathclap{\varrho_c \in \stackrel{\operatorname{c}}{E}\setminus \{\varrho_\gamma\}}}^{} (\varrho_\gamma - \varrho_c) 
	\left[ \frac{1}{|\det g_\gamma |^{1/2}} \partial_{\stackrel{\gamma}{y_i}} g_\gamma^{ij} |\det g_\gamma |^{1/2} \partial_{\stackrel{\gamma}{y_j}} \psi\right.\\
	&\left. - \sum_{q=1}^{r} \frac{1}{|\det g_\gamma |^{1/2}} \partial_{\stackrel{\gamma}{y_i}} g_\gamma^{ij} |\det g_\gamma |^{1/2} \stackrel{\gamma}{\alpha}_{qj} \partial_{t_q} \psi - \sum_{q=1}^{r} g_\gamma^{ij} \stackrel{\gamma}{\alpha}_{qi} \partial_{t_q} \partial_{\stackrel{\gamma}{y_j}} \psi\right.\\
	&\left.+ \sum_{p,q=1}^{r} g_\gamma^{ij} \stackrel{\gamma}{\alpha}_{qi}\stackrel{\gamma}{\alpha}_{pj} \partial_{t_{q}} \partial_{t_p} \psi\right]
	+ \frac{1}{\prod_{\varrho_c \in \stackrel{\operatorname{c}}{E}\setminus\{\varrho_\gamma\}}(\varrho_\gamma - \varrho_c)^{m(\varrho_c)/2-1}} f_\gamma \psi
	\end{split}
\end{equation}
Again: using $\partial_{t_i}\psi=i \omega_i \psi$ gives us the desired result. 

\noindent The last thing to prove is to show that if $\psi$ fulfills \eqref{eq:separated ode} and \eqref{eq:separated pde} for some constants $\tilde{\lambda}_0, \ldots, \tilde{\lambda}_{r+R-1}, \omega_1,\ldots,\omega_r$ then it is also an eigenfunction of $\stackrel{s}{\hat{K}}$ for all real values $s$. To do so, we use $\partial_{t_i}\psi=i \omega_i \psi$ to obtain \eqref{eq:separatedode} and \eqref{eq:separatedpde} from \eqref{eq:separated ode} and \eqref{eq:separated pde}. Then for all non-constant eigenvalues $\varrho_k$ we multiply the corresponding equation of \eqref{eq:separated ode} by 
\begin{equation*}
	\prod_{\varrho_i \in E\setminus \{\varrho_k\}} \left(\frac{s-\varrho_i}{\varrho_k - \varrho_i}\right)^{m(\varrho_i) /2}
\end{equation*}
and for each constant eigenvalue $\varrho_{\gamma}$ we multiply the corresponding equation of \eqref{eq:separatedpde} by 
\begin{equation*}
	\prod_{\varrho_i \in E\setminus \{\varrho_\gamma\}} \left(\frac{s-\varrho_i}{\varrho_k - \varrho_i}\right) \prod_{\varrho_i \in \stackrel{\operatorname{c}}{E}} (s- \varrho_{\gamma})^{m_\gamma /2 -1}
\end{equation*}
and add up the results. Then we use the product rule for differentiation and the lagrange interpolation formula for polynomials and arrive at $\stackrel{s}{\hat{K}} \psi = \sum_{k=1}^{n-1}\lambda_k s^k \psi$. Theorem 8 is proven.

\section{Summary}
We have generalized the integrals that Topalov \cite{topalov2003geodesic} has found for the geodesic flow of {c-compatible} structures $(g,J,A)$ to certain natural hamiltonian systems.
We have proven that these integrals also commute in the quantum sense. The result was then generalized to a class of natural Hamiltonian systems:
in the case where the tensor $A$ is semi-simple we have described all potentials that may be added to the kinetic energy term such that the resulting functions on $T^* \mathcal{M}$ still poisson-commute pairwise and their quantum operators commute pairwise as well. The potentials that are admissible in the quantum problem are the same as in the classic case.
In the case that $A$ is not semi-simple we could present some potentials that may be added to the kinetic energy such that the modified integrals still commute in both the classical and the quantum sense; it is however not clear whether there exist more and this shall be subject to further investigation.
We have tackled the question of the separation of variables for common eigenfunctions of the constructed differential operators in the case where $A$ is semi-simple and all constant eigenvalues are real. If all eigenvalues of $A$ are non-constant we get complete reduction to ordinary differential equations. The case where $A$ has constant non-real eigenvalues or Jordan blocks still needs investigation.
\section{Acknowledgements}
I would like to thank my advisor professor Vladimir Matveev for the suggestion of this topic and the useful discussions.
	
\end{document}